\documentclass[11pt,a4paper]{amsart}

\usepackage{makecell}
\usepackage{mathpazo}
\usepackage[mathscr]{eucal}
\usepackage{enumerate}
\usepackage[psamsfonts]{amssymb}
\usepackage{pdfsync}
\usepackage{graphicx,psfrag}
\usepackage[psamsfonts]{amssymb}
\usepackage{amscd}
\usepackage{amsmath}
\usepackage{amsxtra}
\usepackage{mathrsfs}
\usepackage{stmaryrd}
\usepackage{comment}
\usepackage{enumerate}
\usepackage{hyperref}
\usepackage[capitalize,nameinlink]{cleveref}
\usepackage{tikz-cd}

\usepackage{esint}

\usepackage{subcaption}

\usepackage[width=5.8in, height=8.5in, bottom=1.3in, centering]{geometry}

\usepackage[all]{xy}

\usepackage{mathtools}

\usepackage{mathrsfs}

\usepackage{hyperref}
\hypersetup{ colorlinks, citecolor=green, filecolor=black, linkcolor=blue, urlcolor=blue } 

\theoremstyle{plain}
        \newtheorem{theorem}{Theorem}[section]
        \newtheorem{lemma}[theorem]{Lemma}
        \newtheorem{proposition}[theorem]{Proposition}

        \newtheorem{Theorem}[]{Theorem}
        \theoremstyle{definition}
        \newtheorem{definition}[theorem]{Definition}
        \newtheorem{remark}[theorem]{Remark}
        \newtheorem{example}[theorem]{Example}

\newcommand{\Z}{\mathbb{Z}}

\newcommand{\sfM}{\mathsf{M}}
\newcommand{\sfD}{\mathsf{D}}
\newcommand{\D}{\mathsf{d}}

\newcommand{\sfU}{\mathsf{U}}
\newcommand{\sfOmega}{\mathsf{\Omega}}

\newcommand{\frg}{\mathfrak{g}}

\newcommand{\Hom}{\operatorname{Hom}}
\newcommand{\Der}{\operatorname{Der}}

\newcommand{\ulr}{\sfU(L,R)}
\renewcommand{\hat}{\widehat}

\newcommand{\pbww}{\widetilde{\mathsf{pbw}}}

\renewcommand{\tilde}{\widetilde}

\author{Bjarne Kosmeijer}
\address{University of Amsterdam }
\email{b.a.kosmeijer@uva.nl}
\author{Hessel Posthuma}
\address{University of Amsterdam }
\email{h.b.posthuma@uva.nl}

\title{Hochschild cohomology of Lie--Rinehart algebras}
\begin{document}
\begin{abstract}
We compute the Hochschild cohomology of universal enveloping algebras of Lie--Rinehart algebras
 in terms of the Poisson cohomology of the associated graded quotient algebras. 
Central in our approach are two cochain complexes of ``nonlinear Chevalley--Eilenberg'' cochains whose origins 
lie in Lie--Rinehart modules ``up to homotopy'', one on the Hochschild cochains of the base algebra, another related to the adjoint representation.
The Poincar\'e--Birkhoff--Witt isomorphism is then extended to a certain intertwiner between such modules.
Finally, exploiting the twisted Calabi--Yau structure, we obtain results for the dual Hochschild and cyclic homology.
\end{abstract}
\maketitle

\tableofcontents
\section*{Introduction}
Lie--Rinehart algebras are generalizations of Lie algebras that
first appeared in the work of Rinehart \cite{rinehart} on the Hochschild--Kostant--Rosenberg theorem. 
The main generalization stems from the fact that they are defined over another commutative algebra: a Lie--Rinehart algebra 
is given by a pair $(L,R)$ where $L$ is a Lie algebra that acts on a commutative algebra $R$ (via its {\em anchor}). Furthermore,
$L$ is an $R$-module and there is a Lebniz-type rule showing that the Lie bracket is a derivation over the module structure. (The precise definition will
be recalled below.)

Extreme examples of Lie--Rinehart algebras are given by Lie algebras, --where $R$ is a field and the anchor is zero--, as well as vector fields
on a manifold $M$ where $R=C^\infty(M)$ is given by the ring of smooth functions. More examples with $R=C^\infty(M)$ come from Lie algebroids, 
in fact Lie--Rinehart algebras can be viewed as the algebraic version of Lie algebroids. 

Many aspects of the theory of Lie algebras can be generalized to Lie--Rinehart algebras in a natural way: there is a natural notion of modules, their cohomologies 
and these are all controlled by the {\em universal enveloping algebra} $\sfU(L,R)$. Keeping in mind the example given by vector fields gives a geometric flavor to the 
theory, where modules correspond to flat vector bundles and the universal enveloping algebra is given by the algebra of differential operators. 

In this paper we study the Hochschild cohomology of universal enveloping algebras of Lie--Rinehart pairs. By the Poincar\'e--Birkhoff--Witt (PBW) theorem proved in \cite{rinehart}
the symmetric algebra ${\rm Sym}_R(L)$ inherits a Poisson structure which can be used to define its \emph{Poisson cohomology} $H^\bullet_{\rm Pois}({\rm Sym}_R(L))$, c.f.\ \cite{Huebschmann}.
Our first main result is given as follows:
\begin{Theorem}
\label{thm-m1}
Let $(L,R)$ be a Lie --Rinehart algebra with $R$ smooth and $L$ projective over $R$. There is a natural isomorphism
\[
HH^\bullet(\sfU(L,R),\sfU(L,R))\cong H^\bullet_{\rm Pois}({\rm Sym}_R(L)).
\]
Furthermore, the Poisson cohomology is naturally isomorphic to the Lie--Rinehart cohomology of the symmetric algebra of the adjoint representation:
\[
H^\bullet_{\rm Pois}({\rm Sym}_R(L))\cong H^\bullet(L,{\rm Sym}({\rm ad}(L)).
\]
\end{Theorem}
Remark that the adjoint representation does not exist as an honest Lie--Rinehart module, but only "up to homotopy" as in \cite{ac} constructed using a connection.
The Poisson cohomology complex can be viewed as providing a connection free construction of the representation up to homotopy similar to
the ``nonlinear cochain complex'' for deformation cohomology of Lie algebroids in \cite{crmo}.
The theorem above was first proved in a topological setting in the context of Lie algebroids (i.e., where $R=C^\infty(M)$) in \cite{blom} using Kontsevich' formality theorem and viewing $\sfU(L,R)$ as a quantization of ${\rm Sym}_R(L)$. Furthermore, the isomorphism $HH^\bullet(\sfU(L,R),\sfU(L,R))\cong H^\bullet(L,{\rm Sym}({\rm ad}(L))$ was conjectured in the algebraic setting in \cite{Bouaziz}.

Our approach in this paper is completely different to the approach of Blom and works in the algebraic context as well. It 
builds on the idea to view $\sfU(L,R)$-modules as both an $R$-module and an $L$-module. 
On the level of derived functors computing ${\rm Ext}^\bullet_{\sfU^e}(\sfU,-)$ this leads to writing the functor of $\sfU(L,R)^e$-invariants as a composition of the functors taking
$R^e$-invariants and $L$-invariants as in \cite{KL}, but the main issue in the computation is that after taking projective resolutions, the $R$- and $L$-module structure 
are not compatible on the nose, but only after taking cohomology. We resolve this issue by constructing explicit homotopies controlling this incompatibility and using these
to write down an explicit cochain complex of ``nonlinear Chevalley--Eilenberg cochains''. We embed these homotopies into a notion of a \emph{quasi-modules} over a Lie--
Rinehart algebra, offering a different point of view on the adjoint representation (up to homotopy) and an  explanation for the origins of connection free nonlinear cochain 
models for representations up to homotopy of Lie algebroids. 

In degree one, the theorem above can be viewed as an infinitesimal version of our previous result \cite{kopo} where the deformation cohomology of a 
Lie groupoid is related to the Hochschild cohomology of its convolution algebra. As explained in that paper, these two results can be related explicitly using deformation
quantization and the van Est map.

Our second result concerns the dual Hochschild and cyclic homology of universal enveloping algebras. It was proved recently in \cite{llm} that $\sfU(L,R)$ is twisted 
Calabi--Yau with dualizing bimodule induced from the $(L,R)$-module $Q_L:=\det_R(L)\otimes_R\Omega^{\rm top}_R$ of ``transverse densities''.
With this we prove that the dual homologies are related to the {\em Poisson homology} of ${\rm Sym}_R(L)$, computed from chain complex of algebraic differential forms 
$\Omega^\bullet_{{\rm Sym}_R(L)}$ equipped with the Poisson differential $L_P$:
\begin{Theorem}
\label{thm-m2}
Assume furthermore that $L$ is finitely generated of constant rank $d$ over $R$ and that $R$ satisfies van den Bergh duality in dimension $n$.
Then there are natural isomorphisms
\begin{align*}
H_\bullet(\sfU(L,R))&\cong H_\bullet(\Omega^\bullet_{{\rm Sym}_R(L)},L_P),\\
HC_\bullet(\sfU(L,R))&\cong H_\bullet(\Omega^\bullet_{{\rm Sym}_R(L)}[u],L_P+ud).
\end{align*}
Furthermore, for the Poisson homology we have the isomorphism 
\[
H_\bullet(\Omega^\bullet_{{\rm Sym}_R(L)},L_P)\cong H^\bullet(L,{\rm Sym}({\rm ad}(L,R))\otimes Q_L).
\]
\end{Theorem}
The set up of this paper is as follows: after a brief recollection of some basic notions of Lie--Rinehart algebras in \cref{sec-pr}, we develop the 
theory of quasi-modules in \cref{sec-qm}, and discuss the two main examples: the adjoint representation and the quasi-module of Hochschild cochains.
We then show in \cref{sec:var-pbw} how to extend the Poincar\'e--Birkhoff--Witt (PBW) isomorphism to an intertwiner between these two modules, and use this in 
\cref{sec-hc} to compute the Hochschild cohomology of the universal enveloping algebra and prove \cref{thm-m1}. Finally, in \cref{sec-hm} we combine the 
extended PBW morphism with van den Bergh duality to prove \cref{thm-m2}

\subsection*{Acknowledgement}
The research B.K. is supported by the Netherlands Organization for Scientific Research through project nr. 613.001.021.

\section{Preliminaries}
\label{sec-pr}
\subsection{Lie--Rinehart algebras}
Here we recall the basic definitions in the theory of Lie--Rinehart algebras. All the concepts discussed here can be traced back to the foundational paper \cite{rinehart}. 
We work over a fixed field $\mathbb{K}$ of characteristic $0$.
The central definition of this paper is that of a Lie--Rinehart algebra itself:
\begin{definition}
A {\em Lie--Rinehart algebra} is a pair $(L,R)$ with $R$ a commutative algebra and $L$ a Lie algebra that is also 
equipped with the structure of an $R$-module. Furthermore there is a morphism $\rho:L\to{\rm Der}(R)$ of Lie algebras, called {\em the anchor},
satisfying the {\em Leibniz identity}:
\[
[X, rY]=\rho(X)(r)Y+r[X,Y],\qquad\mbox{for}~r\in R,~X,Y\in L.
\]
\end{definition}
We shall often simply write $X(r)$ for the action $\rho(X)(r)$ of $L$ on $R$ via the anchor.
Clearly, setting $L:={\rm Der}(R)$ for any commutative algebra $R$ produces the first examples of Lie--Rinehart algebras, but below we will
see that there are many more examples.
\begin{definition}
\label{def-lrm}
An $(L,R)$-{\em module} is given by an $R$-module $M$ that is equipped with a Lie algebra module structure given by 
$R$-linear maps $\mathcal{L}:L\to {\rm End}(M),~X\mapsto \mathcal{L}_X$ satisfying
\begin{equation}
\label{eq-leibniz-module}
\mathcal{L}_X(rm):=X(r)m+r\mathcal{L}_X(m)
\end{equation}
\end{definition}
\begin{remark}
\label{rk:con}
For a projective $R$-module $M$ it is always possible to choose a ``$L$-connection'' given by an $R$-linear map $\nabla: L\to {\rm End}(M)$ satisfying \cref{eq-leibniz-module}.
It may fail however to define an $(L,R)$-module structure because the ``curvature'' is non-zero: 
\[
K(X,Y):=[\nabla_X,\nabla_Y]-\nabla_{[X,Y]}\not = 0.
\]
\end{remark}
Remark that the ``trivial module'' is given by $R$ equipped with the $L$-module structure given by the anchor.
Such modules over a Lie--Rinehart algebra $(L,R)$ are in one-to-one correspondence with modules over its {\em universal enveloping algebra} $\sfU(L,R)$
constructed as follows: first form the Lie algebra $\mathfrak{g}:=L\oplus R$ with Lie bracket
\[
[(X,r),(Y,s)]:=([X,Y],X(s)-Y(r)),\qquad\mbox{for}~X,Y\in L,~r,s\in R.
\]
We then form the universal enveloping algebra $\sfU(\mathfrak{g})$ of this Lie algebra and quotient out the two-sided ideal 
generated by $R$-module structure on $L$:
\[
\sfU(L,R):=\sfU(\mathfrak{g})\slash I,\qquad I:=\left<r\cdot X-rX,~r\in R,~X\in L\right>.
\]
In the following, we will sometimes write $\sfU$ for $\sfU(L,R)$, when the underlying Lie-Rinehart algebra is clear from the context. This algebra has the following universal property: let $A$ be an algebra equipped with a pair of morphisms $\phi_R:R\to A$ and $\phi_L:L\to\mathfrak{gl}(A)$
satisfying  
\[
\phi_L(rX)=\phi_R(r)\phi_L(X),\qquad\mbox{for}~r\in R,~X\in L.
\]
Then $\phi_R$ and $\phi_L$ factor over $\sfU(L,R)$ via the canonical maps $i_{R,L}:R,L\to \sfU(L,R)$: there exists a unique morphism $\phi:\sfU(L,R)\to A$ of algebras such that 
the diagrams
\[
\xymatrix{R\ar[rr]^{\phi_R}\ar[d]_{i_R}&&A\\
\sfU(L,R)\ar[urr]_\phi&&}
\qquad 
\xymatrix{L\ar[rr]^{\phi_L}\ar[d]_{i_L}&&A\\
\sfU(L,R)\ar[urr]_\phi&&}
\]
commute. With the correspondence between $(L,R)$-modules and $\sfU(L,R)$ modules the {\em Lie--Rinehart cohomology} with values in a module $M$
is defined as 
\[
H^\bullet_R(L;M):=\mathsf{Ext}^\bullet_{\sfU}(R,M).
\] 
When $L$ is projective as an $R$-module, there is a concrete cochain complex computing this cohomology: the Chevalley--Eilenberg cochains are given by
\begin{equation}
\label{eq-cec}
\sfOmega^\bullet_R(L;M):={\rm Hom}_R(\bigwedge^\bullet_RL,M),
\end{equation}
with differential $\delta_{\mathsf{CE}}:\sfOmega_R^k(L;M)\to \sfOmega_R^{k+1}(L;M)$ defined as
\begin{align}
\label{eq-ced}
(\delta_{\mathsf{CE}}\phi)(X_1,\ldots,X_{k+1})=&\sum_{i=1}^{k+1}(-1)^{i+1}\mathcal{L}_{X_i}(\phi(X_1,\ldots,\widehat{X}_i,\ldots,X_{k+1}))\\
\notag&+\sum_{i<j}(-1)^{i+j}\phi([X_i,X_j],X_1,\ldots,\widehat{X}_i,\ldots,\widehat{X}_j,\ldots,X_{k+1}).
\end{align}

\begin{example}[Lie algebras]
\label{ex-la}
When $R=\mathbb{K}$, a Lie--Rinehart algebra over $k$ is simply given by a Lie algebra $L=\mathfrak{g}$. In this case the anchor $\rho=0$, and all definitions above reduce to the standard concepts in the theory of Lie algebras: a module in the sense of Lie--Rinehart algebras is just a Lie algebra module and we see that the universal enveloping algebra is just given by $\sfU(\mathfrak{g})$. In this case $(\sfOmega^\bullet_\mathfrak{g}(\sfM),\delta_{\mathsf{CE}})$ is the standard Chevalley--Eilenberg cochain complex computing the Lie algebra cohomology $H(\mathfrak{g};\sfM)$.
\end{example}
\begin{example}[Differential operators on affine varieties]
\label{ex-do}
When $R$ is the coordinate ring of a smooth affine variety $\mathsf{X}$, the Lie algebra $L={\rm Der}(R)$ of vector fields 
on $\mathsf{X}$ forms a Lie--Rinehart pair with $\rho={\rm id}$. In this case $\sfU(L,R)$ equals the algebra of differential operators $\mathcal{D}(\mathsf{X})$ on 
$\mathsf{X}$. In the case that $R=\mathbb{K}[x_1,\ldots,x_n]$ this is called the Weyl algebra $A_n$.
A module over this Lie--Rinehart algebra is the same thing as a vector bundle $E$ equipped with a flat connection $\nabla^E$. In this case,
the Chevalley--Eilenberg complex with values in $R$ equals the algebraic de Rham complex $(\sfOmega^\bullet(\mathsf{X}),\D_{\rm dR})$ computing the cohomology of $M$.  
\end{example}
Other examples come from actions of Lie algebras by derivations on a commutative algebra, logarithmic vector fields along a divisor in an affine variety, hyperplane arrangements, foliations etc.
\subsection{The adjoint representation}
\label{sec-ad}
Whereas Lie algebras have a canonical adjoint representation on itself, Lie--Rinehart algebras do not have this luxury. 
This is the case in particular for Lie algebroids, but in \cite{ac} it is shown that Lie algebroids do have a canonical adjoint representation {\em up to homotopy}.
Under the assumption that $L$ is projective as an $R$-module, we can easily generalize this construction to Lie--Rinehart algebras $(L,R)$. 
This assumption implies that we can choose a connection
$\nabla:{\rm Der}(R)\to{\rm End}(L)$. Following \cite{ac} (these authors work in the geometric setting of Lie algebroids), we use this connection to define
the following $L$-connections on $L$ and ${\rm Der}(R)$:
\begin{align*}
\nabla^L_X(Y)&:=\nabla_{\rho(Y)}X+[X,Y],\qquad\qquad X,Y\in L\\
\nabla^{{\rm Der}(R)}_X(D)&:=\rho(\nabla_D(X))+[\rho(X),D],\qquad X\in L,~D\in{\rm Der}(R)
\end{align*}
These connections are compatible in the sense that $\nabla^L\circ \rho=\rho\circ\nabla^{{\rm Der}(R)}$, and are jointly called the {\em basic connection} 
$\nabla^b=(\nabla^L,\nabla^{{\rm Der}(R)})$. 
We introduce the {\em basic curvature} $K^b_\nabla\in \sfOmega^2(L;{\rm Hom}_R({\rm Der}(R),L))$ as
\[
K^b_\nabla(X,Y)(D):=\nabla_D([X,Y])-[\nabla_DX,Y]-[X,\nabla_D(Y)]-\nabla^{{\rm Der}(R)}_{\nabla^L_Y(D)}(X)+\nabla^{{\rm Der}(R)}_{\nabla^L_X(D)}(Y),
\]
for $X,Y\in L$ and $D\in {\rm Der}(R)$. It has the property that $K_\nabla^b\circ \rho$ equals the curvature of $\nabla^L$, whereas $\rho\circ K_\nabla^b$ is the 
curvature of $\nabla^{{\rm Der}(R)}$. We then consider the graded
$R$-module ${\rm ad}(L):=L\oplus{\rm Der}(R)$, where $L$ has degree $0$ and ${\rm Der}(R)$ degree $1$. It is shown in \cite[\S 3.2.]{ac} that the operator
\[
\mathsf{D}:=\rho+\nabla^b+K_\nabla,
\]
defines a differential on $\sfOmega^\bullet_R(L;{\rm ad}(L))$, i.e., $\mathsf{D}^2=0$. This is by definition the adjoint representation up to homotopy, 
and the cohomology of $\mathsf{D}$ is denoted by $H_R^\bullet(L;{\rm ad}(L))$.

Symmetric powers of this adjoint representation are constructed using the tensor product connection, defined on the $R$-module 
${\rm Sym}({\rm ad}(L))=\bigwedge_R{\rm Der}(R)\otimes_R{\rm Sym}_R(L)$.

Although the construction of adjoint representation up to homotopy uses a connection, the cohomology $H_R^\bullet(L;{\rm Sym}({\rm ad}(L)))$ has an intrinsic meaning.
For this we consider the symmetric algebra ${\rm Sym}_R(L)$. This algebra is actually isomorphic to the universal enveloping algebra
if we were to put the Lie-bracket $[~,~]$ and anchor $\rho$ to zero. But for nonzero anchor and bracket, ${\rm Sym}_R(L)$ inherits a Poisson structure given
by
\begin{equation}
\label{eq-epb}
\{r,s\}=0,\quad\{X,r\}=X(r),\quad \{X,Y\}=[X,Y],\qquad r,s\in R,~X,Y\in L,
\end{equation}
and extended to higher degrees by the derivation property of the Poisson bracket. For any Poisson algebra, there exists a Poisson cohomology theory defined in 
\cite{Huebschmann}, generalizing the geometric theory in \cite{Lichnerowicz}. 
In our case of ${\rm Sym}_R(L)$ it is defined using the cochain complex with 
\[
\mathfrak{X}^\bullet_{{\rm Sym}_R(L)}:=\bigwedge^\bullet_R({\rm Der}({\rm Sym}_R(L)),
\]
with differential $\delta_P:\mathfrak{X}^k_{{\rm Sym}_R(L)}\to \mathfrak{X}^{k+1}_{{\rm Sym}_R(L)}$ defined by
\begin{align}
\notag
\delta_P(\mathsf{D})(\gamma_1,\ldots,\gamma_{k+1}):=&\sum_{i=1}^{k+1}(-1)^i\{\gamma_i,\mathsf{D}(\gamma_1,\ldots,\widehat{\gamma}_i,\ldots,\gamma_{k+1})\}\\
\label{eq-pdif}
&+\sum_{i<j}(-1)^{i+j}\mathsf{D}(\{\gamma_i,\gamma_j\},\gamma_1,\ldots,\widehat{\gamma}_i,\ldots,\widehat{\gamma}_j,\ldots,\gamma_{k+1}).
\end{align}
Consider now the map 
\begin{equation}
\label{eq-nlm}
\mathfrak{X}^p_{{\rm Sym}_R(L)}\to \bigoplus_{i=0}^p {\rm Hom}\left(\bigwedge^{p-i} L,\bigwedge^i_R {\rm Der}(R)\otimes_R{\rm Sym}_R(L)\right),
\end{equation}
given by mapping $\mathsf{D}\mapsto (\phi^\mathsf{D}_0,\ldots \phi^\mathsf{D}_p)$, where 
\[
\phi^\mathsf{D}_i(X_1,\ldots,X_{p-i},r_1,\ldots r_i):=\mathsf{D}(X_1,\ldots,X_{p-i},r_1,\ldots r_i).
\]
By the derivation property of $\mathsf{D}$ these maps satisfy
\begin{align}
\notag
\phi^\mathsf{D}_i(X_1,\ldots,X_{p-i-1},rX_{p-i},r_1,\ldots,r_i)=&r\phi^\mathsf{D}_i(X_1,\ldots,X_{p-i-1},X_{p-i},r_1,\ldots,r_i)\\
\label{eq-nlpm}
&+\phi^\mathsf{D}_{i+1}(X_1,\ldots,X_{p-i-1},r,r_1,\ldots,r_i)X_{p-i}.
\end{align}
It is exactly this type on nonlinearity that we shall encounter also in the computation of the Hochschild cohomology below. Notice that because 
${\rm Sym}_R(L)$ is generated by $R$ and $L$, one easily observes that the map in \cref{eq-nlm} is an isomorphism onto the subspace of tuples of maps
satisfying the conditions in \cref{eq-nlpm}. One can write out the Poisson differential in \cref{eq-pdif} in terms of the cochains $(\phi_0,\ldots \phi_p)$:
for the first one we simply find $\delta_{\mathsf{CE}}\phi_0$, the higher ones can be recovered from this one using the condition in \cref{eq-nlpm}. This observation leads to a complex that we will call $\sfOmega_{\rm nl}({\rm Sym}({\rm ad}(L)))$, which will be introduced from another perspective in \cref{section:epc} and which by construction is isomorphic to $\mathfrak{X}^\bullet_{{\rm Sym}_R(L)}$. We remark that it is akin to the model of Abad and Crainic \cite[Example 4.5]{ac} for the anti-symmetric powers of the adjoint representation up to homotopy in the geometric case.

To connect this complex with the symmetric powers of the adjoint representation in the algebraic setting, first notice the following
\begin{lemma}
\label{lemma-ses}
Let $L$ be projective as a module over $R$.
There exists a short exact sequence of ${\rm Sym}_R(L)$-modules
\[
0\longrightarrow {\rm Sym}_R(L)\otimes_R L^\vee\longrightarrow {\rm Der}({\rm Sym}_R(L))\longrightarrow {\rm Der}(R,{\rm Sym}_R(L))\longrightarrow 0,
\]
and choosing splitting is the same as choosing a connection $\nabla$ on $L$.
\end{lemma}
\begin{proof}
Let $\pi:{\rm Sym}_R(L)\to L$ be the projection onto the degree $1$ part of the symmetric algebra. 
Given a connection $\nabla$, construct a map $\sigma:{\rm Der}(R,{\rm Sym}_R(L))\to {\rm Der}({\rm Sym}_R(L))$ as follows:
Given $D\in {\rm Der}(R,{\rm Sym}_R(L))$, let $D_0:=\pi\circ D$. Then we define
\[
\sigma(D)(r):=D(r),\quad \sigma(D)(X):=\nabla_{D_0}(X),
\]
and we uniquely extend this to higher degrees using the derivation property. This defines a splitting.

Conversely, given a splitting $\sigma:{\rm Der}(R,{\rm Sym}_R(L))\to {\rm Der}({\rm Sym}_R(L))$, define
\[
\nabla_D(X):=\pi(\sigma(D)(X)),\qquad D\in{\rm Der}(R),~X\in L.
\]
It is easy to check that this defines a connection. 
\end{proof}
We can now connect with the adjoint representation up to homotopy. Let $(c_0,\ldots,c_k)\in\sfOmega_R(L;{\rm Sym}({\rm ad}(L)))$ be a cochain 
of degree $p$, with $c_i\in{\rm Hom}_R(\bigwedge^{p-i}_RL,\bigwedge^i_R{\rm Der}(R)\otimes_R{\rm Sym}_R(L))$. We then define recursively 
\begin{align*}
\phi_i(X_1,\ldots,X_{p-i},r_1,\ldots, r_i):=&c_i(X_1,\ldots,X_{p-i},r_1,\ldots, r_i)\\
&-\sum_{j=1}^{p-i}(-1)^j\sigma(c_{i+1}(X_1,\ldots,\hat{X}_j,\ldots,X_{p-i},-,r_1,\ldots,r_i))(X_j)
\end{align*}
This defines an isomorphism of cochain complexes
\[
\sfOmega_{\rm nl}(L;{\rm Sym}({\rm ad}(L))) \cong \sfOmega_R(L;{\rm Sym}({\rm ad}(L))),
\]
and this proves the second statement in \cref{thm-m1}.

\subsection{Hochschild cohomology}
Let $(L,R)$ be a Lie--Rinehart algebra with universal enveloping algebra $\sfU(L,R)$. In the following we shall sometimes simply write $\sfU$ for this algebra
when the underlying Lie--Rinehart pair is clear. Let $\mathsf{M}$ be a bimodule over $\sfU$.
The Hochschild cohomology groups $HH^\bullet(\sfU(L,R),\mathsf{M}):=\mathsf{Ext}_{\sfU^e}^\bullet(\sfU,\mathsf{M})$ are  the derived functors of
\[
{\rm Hom}_{\sfU^e}(\sfU,-):\mathsf{Mod}_{\sfU^e}\to \mathsf{Vect}_k,
\]
applied to $\mathsf{M}$. 
Using the bar-resolution of $\sfU$, we obtain the standard Hochschild complex given by $C^k(\sfU,\mathsf{M}):={\rm Hom}(\sfU^{\otimes k},\mathsf{M})$ with
differential 
\begin{align}
\label{eq-hd}
b\phi(\sfD_1,&\ldots,\sfD_{k+1}):=\sfD_1\phi(\sfD_2,\ldots,\sfD_{k+1})\\
\notag &+\sum_{i=1}^k(-1)^i\phi(\sfD_1,\ldots, \sfD_i\sfD_{i+1},\ldots,\sfD_k)+(-1)^{k+1}\phi(\sfD_1,\ldots,\sfD_k)\sfD_{k+1}.
\end{align}
We are mainly interested in the case $\mathsf{M}=\sfU$. In this case $C^\bullet(\sfU,\sfU)$ 
has a rich algebraic structure and becomes a differential graded algebra by means of the cup-product:
\begin{equation}
\label{eq-cup}
(\phi_1\cup\phi_2)(\sfD_1,\ldots\sfD_{k_1+k_2}):=\phi_1(\sfD_1,\ldots,\sfD_{k_1})\phi_2(\sfD_{k_1+1},\ldots,\sfD_{k_1+k_2}),\quad\phi_i\in C^{k_i}(\sfU,\sfU).
\end{equation}
\begin{remark}
The universal enveloping algebra $\sfU(L,R)$ has the rich algebraic structure of a \emph{left Hopf algebroid} \cite{kp,CRvdB}. 
Associated to this there exists a ``Hopf--Hochschild and cyclic'' (co)homology, defined on a certain subcomplex of the Hochschild cochain complex that we consider here and
computed in {\em loc.cit.}. In this paper we will not consider the Hopf-structure and consider only the Hochschild theory of $\sfU(L,R)$ as an algebra.
\end{remark}

\section{Quasi-Modules}
\label{sec-qm}
Let $(L,R)$ be a Lie--Rinehart algebra. In this section we analyse a certain type of $L$-modules that we encounter in the computation of the Hochschild 
cohomology of $\sfU(L,R)$ below. Recall \cref{def-lrm} of an $(L,R)$-module. The $L$-modules that we consider here have an additional $R$-module structure,
but the two module structures are not quite compatible with each other to form an $(L,R)$-module.
\begin{definition}
\label{def-qm}
An $(L,R)$-{\em quasi-module} is a a cochain complex $\sfM=\bigoplus_{k\in\Z}\mathsf{M}^k$ with differential 
$\D:\sfM^k\to \sfM^{k+1}$, equipped with an $R$-module structure, denoted by $m\mapsto r\cdot m,~r\in R,~m\in \sfM$,
and a $L$-module structure $\mathcal{L}_X:\sfM\to \sfM,~X\in L$, that are compatible with the differential:
\begin{align*}
r\circ \D&=\D\circ r,\qquad \mbox{for all}~r\in R,\\
\mathcal{L}_X\circ \D&=\D\circ \mathcal{L}_X,\qquad \mbox{for all}~X\in L.
\end{align*}
The action of $R$ and $L$ are compatible via the Leibniz rule:
\[
\mathcal{L}_X\circ r=r\circ \mathcal{L}_X+X(r),
\]
and there exist a map $h:R\otimes L\to{\rm End}_{-1}(\sfM)$,~$(r,X)\mapsto h_{r,X}:\sfM^\bullet\to \sfM^{\bullet-1}$ satisfying
\begin{equation}
\label{def-htp}
\mathcal{L}_{rX}=r\circ \mathcal{L}_X+h_{r,X}\circ \D+\D\circ h_{r,X}.
\end{equation}
Finally, the homotopy operators should satisfy
\begin{align}
\label{prop-htp}
\mathcal{L}_Y\circ h_{r,X}&=h_{r,X}\circ\mathcal{L}_Y+h_{Y(r),X}+h_{r,[Y,X]},
%r_1\circ h_{r_2,X}&=h_{r_2,X}\circ d.
\end{align}
\end{definition}
%\begin{remark}
%Unaware to the authors these objects appeared in \cite{Bouaziz} under the name of derived representations. The property defined by \eqref{prop-htp} corresponds to the contraction map in \cite{Bouaziz} being a morphism of pre-representations.
%\end{remark}
\begin{remark}
\label{rk-htp}
When $h_{r,X}=0$ for all $r\in R,~X\in L$, this is just the definition of a $(L,R)$-module structure on the cochain complex $(\sfM,\D)$.
It follows immediately that the cohomology $H^\bullet(\sfM)$ carries the structure of a graded $(L,R)$-module.
\end{remark}
\begin{remark}
In the literature on Lie algebroids, the operator $\mathcal{L}$ above satisfying \cref{def-htp} has appeared before under the name of 
{\em flat connections up to homotopy}, c.f. \ \cite[\S 3]{cf-sc}. 
The notion of quasi-module as defined above appeared in \cite{Bouaziz} under the name of derived representations, 
the property defined by \eqref{prop-htp} corresponds to the contraction map in \cite{Bouaziz} being a morphism of pre-representations.
Quasi-modules differ from the by now more standard notion of  ``representations up to homotopy'', c.f.\ 
\cite{ac,Block}, which typically depends on the choice of a connection, and it remains to investigate the precise relation between the two. 
It would be interesting to see if and if so, how, these quasi-modules fit
into the theory of \cite{lv}, where a very general representation theory of Lie--Rinehart algebras up to homotopy is outlined.
\end{remark}

%\begin{lemma}
%The homotopy operators satisfy
%\begin{align*}
%\left\{\D,\nabla_Y\circ h_{r,X}-h_{r,X}\circ\nabla_Y-h_{Y(r),X}+h_{r,[Y,X]}\right\}&=0,\\
%\left\{r_1\circ h_{r_2,X}+h_{r_1,r_2X}-h_{r_1r_2,X},\mathsf{d}\right\}&=0.
%\end{align*}
%\end{lemma}
For a quasi module $(\sfM,\D)$,
it is possible for such a module to write down the Chevalley--Eilenberg cochains as in \cref{eq-cec}, but
the incompatibility between the $R$ and $L$-module structure shows that these are not closed under the differential \cref{eq-ced}. To counteract we introduce ``non-linar Chevalley-Eilenberg cochains'' in two steps.

For a given quasi-module $(\sfM,d^\sfM$), we first look at the Chevalley-Eilenberg chains which are linear only with respect to the ground field, i.e. we look at $\sfOmega^{\bullet,\bullet}(L,\sfM)$ with
\[\sfOmega^{p,q}(L,\sfM)={\rm Hom}(\bigwedge^p L,\sfM^q).\]

On this grid, we have two differentials coming from the structure of $(M,d^M)$ as a differentially graded $L$-representation, the Chevalley-Eilenberg differential
\[\delta_\mathsf{CE}: \sfOmega^{p,q}(L,\sfM)\to\sfOmega^{p+1,q}(L,\sfM)\]
given by formula in \cref{eq-ced} and the differential
\[ d^\sfM: \sfOmega^{p,q}(L,\sfM)\to\sfOmega^{p,q+1}(L,\sfM)\]
given by post-composition with the differential of $(\sfM,d^\sfM)$. Since the $L$-action commutes with $d^\sfM$, these differentials commute. We denote by $(\sfOmega^\bullet(L,\sfM),\delta_\mathsf{CE}\pm d^\sfM)$ the associated total complex.

Within this framework we introduce non-linear cochains in the following way:
%\begin{definition}
%Let $(\M,\D)$ be a quasi-$(L,R)$-module. The space of {\em nonlinear Chevalley-Eilenberg cochains} $\Omega^k_{\rm nl}(L,M)$ in degree $k$ is given by
%by those $\phi\in{\rm Hom}\left(\bigwedge^k L, M\right)$ for which there exists a {\em symbol} $\sigma_\phi\in{\rm Hom}\left(\bigwedge^{k-1} L, M\right)$
%such that
%\[
%\phi(rX_1,X_2,\ldots,X_k)=r\phi(X_1,X_2,\ldots,X_k)+dh_{r,X_1}(\sigma_\phi(X_2,\ldots,X_k))+h_{r,X_1}(d\sigma_\phi(X_2,\ldots,X_k)).
%\]
%\end{definition}
%\begin{lemma}
%For a nonlinear cochain $\phi$ with symbol $\sigma_\phi$, we have
%\begin{align*}
%\sigma_{d_{\mathsf{CE}}\phi}&=d_{\mathsf{CE}}\sigma_\phi+\phi,\\
%\sigma_{d^M\phi}&=d^M\sigma_\phi.
%\end{align*}
%\end{lemma}
%\begin{proof}
%To compute the symbol of $d_{\mathsf{CE}}\phi$ we write out the definition
%\begin{align*}
%d_{\mathsf{CE}}\phi(rX_1,\ldots,X_{k+1})=&\nabla_{rX_1}\phi(X_2,\ldots,X_{k+1})\\
%&+\sum_{i=2}^{k+1}(-1)^{i+1}\nabla_{X_i}\phi(rX_1,X_2,\ldots,\hat{X}_i,\ldots,X_{k+1})\\
%&+\sum_{i=2}^{k+1}(-1)^i\phi([rX_1,X_i],X_2,\ldots,\hat{X}_i,\ldots, X_{k+1})\\
%&+\sum_{1<i<j}(-1)^{i+j}\phi([X_i,X_j],rX_1,X_2,\ldots,\hat{X}_i,\ldots,\hat{X}_j,\ldots,X_{k+1}).
%\end{align*}
%Writing out the first line with \cref{def-htp} leads to the second term of the first equation for the symbol in the Lemma.
%The remaining three lines produce $d_{\mathsf{CE}}\sigma_\phi$ for the symbol, if one uses \cref{prop-htp} in the second and third line
%to commute $\nabla_{X_i}$ with $h_{r,X_1}$. The second equality is trivially verified.
%\end{proof}
\begin{definition}
Let $(\sfM,\D^\sfM)$ be a quasi $(L,R)$-module. A {\em nonlinear Chevalley--Eilenberg cochain} in degree $k$ is given by
a sequence $(\phi_0,\ldots,\phi_k)$, with $\phi_i\in{\rm Hom}(\bigwedge^{k-i}L,\sfM^i)$ satisfying
\begin{equation*}
\label{eq-nlcec}
\phi_i(X_1,X_2,\ldots,rX_{k-i})-r\phi_i(X_1,\ldots,X_{k-i})=h_{r,X_{k-i}}(\phi_{i+1}(X_1,\ldots,X_{k-i-1}))
%\quad\mbox{mod}~{\rm Im}(\D^\M).
\end{equation*}
We write $\sfOmega^k_{\rm nl}(L,\sfM)$ for the space of such nonlinear cochains. We thus have inclusions
\[
\sfOmega^\bullet_R(L,\sfM)\subset \sfOmega^\bullet_{\rm nl}(L,\sfM)\subset \sfOmega^\bullet(L,\sfM),
\]
of $R$-linear cochains into non-linear cochains into Chevalley-Eilenberg cochains for the Lie algebra $L$.
As above, we write $\delta_{\rm CE}$ for the Chevalley--Eilenberg differential 
given in \cref{eq-ced}. Since $\sfM$ is not an honest $(L,R)$-module, this differential does not restrict to $\sfOmega^\bullet_R(L;\sfM)$. 
For the nonlinear cochains however, we have:
\end{definition}
\begin{lemma}
The space $\sfOmega^\bullet_{\rm nl}(L;\sfM)=\bigoplus_{k} \sfOmega^k_{\rm nl}(L;\sfM)$ is a subcomplex of $(\sfOmega^\bullet(L,\sfM),\delta_{\mathsf{CE}}\pm d^{\sfM})$.
\end{lemma}
\begin{proof}
This follows from a straightforward computation using the definition \cref{eq-ced} of the Chevalley--Eilenberg differential together with the identities \cref{def-htp,prop-htp}
stated in \cref{def-qm}.
\end{proof}
%It follows from the Lemma that we obtain a double complex $\sfOmega^\bullet_{\rm nl}(L;\mathsf{M}^\bullet)$ with 
%differentials $\mathsf{d}^\mathsf{M}$ and $\delta_{\mathsf{CE}}$. 
The cohomology of this cochain complex
is what we call the {\em Lie--Rinehart cohomology with values in $\mathsf{M}$}, denoted $H_{R}(L;\mathsf{M})$.
%It obviously reduces to the usual Lie--Rinehart cohomology when $\mathsf{M}$ is an honest $(L,R)$-module, c.f. \cref{rk-htp}
\begin{proposition}
Let $(\sfM,\D)$ be a differential graded $(L,R)$-module, i.e., a quasi $(L,R)$-module with $h=0$. Then the complex 
$(\sfOmega^\bullet_{\rm nl}(L,\sfM),\delta_{\mathsf{CE}}\pm \D)$ is isomorphic to the $R$-linear Chevalley-Eilenberg complex $\sfOmega^\bullet_R(L;\sfM)$.
\end{proposition}
\begin{proof}
This is immediate, since plugging $h=0$ into the non-linearity equation \eqref{eq-nlcec} turns it into the equation specifying that the cochains $\phi_i$ are $R$-linear.
\end{proof}
\begin{proposition}
\label{prop-invqi}
Let $q:(\sfM,\D^\sfM)\to (\mathsf{N},\D^{\mathsf N})$ be a quasi-isomorphism between two quasi $(L,R)$-modules that is compatible with the module structures, i.e.,
\begin{align*}
r\circ q&=q\circ r,\\
\nabla_X\circ q&=q\circ \nabla_X,\\
h_{r,X}\circ q&=q\circ h_{r,X}.
\end{align*}
Then $q$ induces a quasi-ismorphism $\tilde{q}:\sfOmega^\bullet_{\rm nl}(L,\sfM)\to \sfOmega^\bullet_{\rm nl}(L,\mathsf{N})$ of cochain complexes.
\end{proposition}
\begin{proof}
Obviously, the map ${\rm Hom}(\bigwedge^i L;\sfM)\to {\rm Hom}(\bigwedge^i L;\mathsf{N})$ given by composition with $q$ restricts, by the properties stated in the Proposition, to a map 
$\tilde{q}:\mathsf{\Omega}^\bullet_{{\rm nl}}(L,\sfM)\to \mathsf{\Omega}^\bullet_{{\rm nl}}(L,\mathsf{N})$. 
Clearly, this map is a morphism of cochain complexes, and compatible with the column-filtration introduced in the proof of the previous proposition.
Therefore the map $\tilde{q}$ induces a morphism of spectral sequences and on the first page
\[
\tilde{q}_1:\mathsf{\Omega}^{\bullet}_R(L;H^{\bullet}(\sfM))\to \mathsf{\Omega}^{\bullet}_R(L;H^{\bullet}(\mathsf{N})),
\] 
is an isomorphism by assumption. It therefore follows that $\tilde{q}$ induces an isomorphism on cohomology.
\end{proof}
\begin{remark}[Tensor products] 
Given two quasi $(L,R)$ modules $(\sfM,\D^\sfM)$ and  $(\mathsf{N},\D^{\mathsf N})$, it is straightforward to equip their tensor product $\sfM\otimes_R\mathsf{N}$ 
with a quasi $(L,R)$ module structure using the product connection $\mathcal{L}_X^{\sfM\otimes_R\mathsf{N}}=\mathcal{L}_X^{\sfM}\otimes 1+1\otimes\mathcal{L}_X^\mathsf{N}$
and homotopies $h_{r,X}^{\sfM\otimes_R\mathsf{N}}=h_{r,X}^{\sfM}\otimes 1+1\otimes h_{r,X}^\mathsf{N}$. This defines a tensor product on the category of quasi-modules. We can also consider (anti-)symmetric tensor products.
\end{remark}
\subsection{Example: the adjoint representation}
\label{section:epc}
Given a Lie--Rinehart algebra $(L,R)$, we consider as before the two-term complex ${\rm ad}(L):=L\oplus{\rm Der}(R)$ with $L$ has degree zero and the differential 
is given by the anchor. We then consider the operator $\mathcal{L}$ simply given by
\[
\mathcal{L}_X(Y):=[X,Y],\qquad \mathcal{L}_X(D):=[\rho(X),D],\qquad\mbox{for}~X,Y\in L,~D\in {\rm Der}(R).
\] 
Together with the homotopy
\[
h_{r,X}(D):=D(r)X,
\]
this defines a quasi-module structure defining the adjoint representation, c.f. \cite[\S 3.3]{cf-sc}. We now consider the symmetric powers of this module. 
This is given by  the complex 
\begin{align*}
{\rm Sym}({\rm ad}(L)):=\bigwedge^\bullet_R {\rm Der}(R)\otimes_R {\rm Sym}_R(L),
\end{align*}
with differential given by
\[
\delta (D_1\wedge\ldots\wedge D_p\otimes X_1\cdots X_k):=\sum_{i=1}^{k}\rho(X_i)\wedge D_1\wedge\ldots\wedge D_p\otimes X_1\cdots\hat{X}_i\cdots X_p.
\]
We give ${\rm Sym}({\rm ad}(L))$  the canonical $R$-module structure as well as an $L$-module structure by
\begin{align*}
\mathcal{L}_X(D_1\wedge\ldots\wedge D_p\otimes X_1\cdots X_k):=&\sum_{i=1}^pD_1\wedge\ldots\wedge [\rho(X),D_i]\wedge\ldots\wedge D_p\otimes X_1\cdots X_k\\
&+D_1\wedge\ldots\wedge D_p\otimes \{X,X_1\cdots X_k\}.
\end{align*}
These two module structures are not compatible, but if we define the homotopy $h_{r,X}:{\rm Sym}^k({\rm ad}(L))\to {\rm Sym}^{k}({\rm ad}(L))$ by
\begin{equation}
\label{eq-h}
h_{r,X}(D_1\wedge\ldots D_p\otimes X_1\cdots X_k):=\sum_{i=1}^p(-1)^{i+1}D_i(r)D_1\wedge\ldots\wedge\hat{D}_i\wedge\ldots\wedge D_p\otimes X_1\cdots X_kX,
\end{equation}
we have:
\begin{proposition}
The triple $({\rm Sym}({\rm ad}(L)),-\delta,h)$ defines an $(L,R)$-module up to homotopy and its complex of nonlinear Chevalley--Eilenberg cochains 
$\sfOmega^\bullet_{\rm nl}(L,{\rm Sym}({\rm ad}(L)))$ is precisely that of Poisson cohomology in \cref{sec-ad}.
\end{proposition}
\begin{proof}
It is easy to check that $\delta$ is compatible with the $R$-module structure and that $\delta^2=0$. To check \cref{def-htp} we remark that 
$r\circ\mathcal{L}_X-\mathcal{L}_{rX}$ acting on $D_1\wedge\ldots \wedge D_p\otimes X_1\cdots X_q\in \bigwedge^p_R {\rm Der}(R)\otimes_R {\rm Sym}^q_R(L)$
is given by
\begin{align*}
\sum_{i=1}^p(-1)^{i+1}D_i(r)\rho(X)\wedge D_1\wedge\ldots\wedge\hat{D}_i\wedge\ldots\wedge D_p&\otimes X_1\cdots X_q\\
+\sum_{j=1}^qX_j(r)D_1\wedge\ldots \wedge D_p&\otimes X_1\cdots\hat{X}_j\cdots X_qX.
\end{align*}
It is easy to check that this expression is exactly $h_{r,X}\circ\delta+\delta\circ h_{r,X}$ acting on $D_1\wedge\ldots \wedge D_p\otimes X_1\cdots X_q$.
Finally, with the homotopies given by \cref{eq-h}, it is clear that the nonlinear Chevalley--Eilenberg cochains are precisely as in \cref{eq-nlpm}.
\end{proof}
\subsection{Example: Hochschild cochains}
\label{sec-h}
Our second example of a quasi $(L,R)$-module comes from Hochschild cohomology.
Let $\sfM$ be a $\sfU(L,R)$-bimodule. By general principles, c.f. \cite[Cor. 2.2]{KL} there exists an $(L,R)$ module structure on the Hochschild cohomology 
$H^\bullet(R,\sfM):={\rm Ext}^\bullet_{R^e}(R,\sfM)$, but the issue is that
it cannot be constructed on the level of the standard Hochschild cochain complex $C^\bullet(R,\sfM)$ after choosing the Bar-resolution 
$\mathsf{Bar}_\bullet(R)$ of $R$ to compute the Ext-groups. In \cite[\S2]{KL} it is explained how to construct an $L$-module structure on $C^\bullet(R,\sfM)$,
but this module structure may fail to be compatible with the $R$-module structure to form a cochain complex of $(L,R)$-modules.
Explicitly, we write $b:C^\bullet(R,\sfM)\to C^{\bullet-1}(R,\sfM)$ for the Hochschild differential given by the formula \cref{eq-hd} restricted from $\sfU$ to $R$ and define for $\varphi\in C^q(R,\sfM)$
\begin{subequations}
\begin{align}
\label{ar}
(r\cdot\varphi)(r_1,\ldots,r_q)&:=r\varphi(r_1,\ldots,r_q),\qquad r,r_1,\ldots,r_q\in R,\\
\label{acl}
\mathcal{L}_X(\varphi)(r_1,\ldots,r_q)&:=\sum_{i=1}^q\varphi(r_1,\ldots,X(r_i),\ldots, r_q)+[X,\varphi(r_1,\ldots,r_q)],\qquad X\in L.
\end{align}
\end{subequations}
It is easy to see that with these two module structures, the Leibniz rule \cref{eq-leibniz-module} 
is satisfied, but $\mathcal{L}_{rX}\not = r\circ\mathcal{L}_X$ for $r\in R$ and $X\in L$. However if we define the operators $h_{r,X}:C^\bullet(R,\sfM)\to C^{\bullet-1}(R,\sfM)$ by
\begin{align}
\label{eq-hqa}
h_{r,X}(\varphi)(r_1,\ldots,r_{q-1})=&\sum_{i<j}^{q-1}\varphi(r_1,\ldots,r_{i-1},r,r_i,\ldots ,X(r_j),\ldots,r_{q-1})\\
\notag
&+\sum_{i=1}^q\varphi(r_1,\ldots,r_{i-1},r,r_{i},\ldots,r_q)X,
\end{align}
we have the following:

\begin{proposition}
The operators in \cref{ar,acl,eq-hqa} give the Hochschild cochain complex 
$(C^\bullet(R,\sfM),b)$ the structure of a quasi $(L,R)$-module as in \cref{def-qm}. 
\end{proposition}
\begin{proof}
It is easy to verify that the identities
\begin{align*}
[\mathcal{L}_X,b]&=0,\\
[\mathcal{L}_X,\mathcal{L}_Y]&=\mathcal{L}_{[X,Y]},\\
\mathcal{L}_X\circ r&=r\circ \mathcal{L}_X+X(r)\circ\mathcal{L}_X,\\
%\nabla_{rX}&=r\circ \nabla_X+h_{r,X}\circ b+b\circ h_{r,X},
\end{align*}
hold true. Next we need to show that
\[
\mathcal{L}_{rX}=r\circ \mathcal{L}_X+h_{r,X}\circ b+b\circ h_{r,X}.
\]
This can be proved by an explicit computation, but to give some insight into the structure
let us divide the derivation of the homotopies into two parts according to the two terms in \cref{acl}. For the first part consider the bar resolution $\mathsf{Bar}_\bullet(R):=R^{\otimes(\bullet+2)}$ and let $L$ act by
\[
\mathcal{L}^1_X(r_0\otimes\ldots\otimes r_{q+1}):=\sum_{i=0}^{q+1}r_0\otimes\ldots\otimes X(r_i)\otimes\ldots\otimes r_{q+1}.
\] 
Together with the $R$ action given by
\[
r\cdot (r_0\otimes\ldots\otimes r_{q+1}):= (rr_0)\otimes r_1\otimes \ldots\otimes r_{q+1},
\]
this induces the action \cref{ar,acl}, except for the last term in $\nabla_X$. We now claim that
\[
\mathcal{L}^1_{rX}-r\circ\mathcal{L}^1_X=h^1_{r,x}\circ \beta+\beta\circ h^1_{r,X},
\]
with 
\begin{align*}
h^1_{r,X}(r_0\otimes\ldots\otimes r_{q+1})=&\sum_{0<i<j\leq q+1}^{q+1}r_0\otimes r_1\otimes\ldots \otimes r_{i-1}\otimes r\otimes r_i\otimes\ldots\otimes X(r_j)\otimes\ldots\otimes r_{q+1}\\
\end{align*}
For the last term in \cref{acl}, we consider the action
\[
\mathcal{L}^2_X(\varphi)(r_1,\ldots,r_q):=[X,\varphi(r_1,\ldots,r_q)],\qquad\mbox{for}~\varphi\in C^q(R,\sfM).
\]
For this action we have as well that
\[
\mathcal{L}^2_{rX}-r\circ\mathcal{L}^2_X=h^2_{r,x}\circ \beta+\beta\circ h^2_{r,X},
\]
with
\[
h^2_{r,X}(\varphi)(r_1,\ldots,r_{q-1}):=\sum_{i=1}^q\varphi(r_1,\ldots,r_{i-1},r,r_{i},\ldots,r_q)X.
\]
Together, $h^1_{r,X}$ and $h^2_{r,X}$ give the homotopy stated in \cref{eq-hqa} under the 
isomorphism ${\rm Hom}_{R^e}(\mathsf{Bar}_q(R),\sfM)\cong C^q(R,\sfM)$ given by setting
$\varphi(r_1,\ldots,r_q):=\tilde{\varphi}(1,r_1,\ldots,r_q,1)$ for $\tilde{\varphi}\in {\rm Hom}_{R^e}(\mathsf{Bar}_q(R),\sfM)$. (Notice that some of the terms in $h^1_{r,X}$ drop out because $X(1)=0$.)

The final identity 
\[
\mathcal{L}_Y\circ h_{r,X}=h_{r,X}\circ\mathcal{L}_Y+h_{Y(r),X}+h_{r,[Y,X]}
\]
is easily verified from the expression in \cref{eq-hqa}. This finishes the proof.
\end{proof}
\section{Variations on PBW}
\label{sec:var-pbw}
The fundamental tool in our computation of Hochschild cohomology of the universal 
enveloping algebra $\sfU(L,R)$ is the Poincar\'e--Birkhoff--Witt isomorphism with the 
symmetric algebra ${\rm Sym}_R(L)$, going back to the fundamental paper of Rinehart 
\cite{rinehart}. This isomorphism is an extention of the well-known PBW map for Lie algebra,
but has the property that in the setting of Lie--Rinehart algebras it is {\em not} a morphism of $L$-modules. In fact, this is closely related to the absence of an adjoint representation as we discussed in \cref{sec-ad}. 

In this section we shall extend the PBW isomorphism for Lie--Rinehart algebras to a quasi-isomorphism of quasi-modules, relating the two main examples of \cref{section:epc} and \cref{sec-h}. In this setting we then prove that this extended PBW-map is a morphism of 
$L$-modules up to a homotopy that allows a further extension to the Chevalley--Eilenberg complex.

We introduce the following notation, we shall write $\vec{D}\otimes\vec{X}$ for a generic element
in $\bigwedge^p_R{\rm Der}(R)\otimes_R{\rm Sym}_R^q(L)$, with 
\[
\vec{D}:=D_1\wedge\ldots\wedge D_p,\qquad \vec{X}:=X_1\cdots X_q,\qquad D_i\in{\rm Der}(R),~X_j\in L.
\]
We then write 
\[
\vec{D}_{(i)}:=D_1\wedge\ldots\wedge\hat{D}_i\wedge\ldots\wedge D_p,\qquad \hat{X}_{(i)}:=X_1\cdots \vec{X}_i\cdots\vec{X}_q,
\]
for the elements with the $i$'th component deleted. Furthermore, $\vec{X}_{(i,j)}$ denotes the same element with the $i$'th and $j$'th component omitted, etc. 

\subsection{The Poincar\'e--Birkhoff--Witt isomorphism for Lie--Rinehart algebras}
\label{section:pbw}
A fundamental tool to understand the structure of the universal enveloping algebra is given by the Poincar\'e--Brikhoff--Witt isomorphism as first proved by Rinehart in \cite{rinehart}: Introduce a filtration given by
\begin{equation}
\label{eq-pbwf}
F_0(U(L,R))=R,\quad F_k(U(L,R)):={\rm span}_\mathbb{K}\{i_L(L)^l,~l\leq k\},~k=1,2,\ldots.
\end{equation}
With this, $U(L,R)$ becomes a filtered algebra; $F_0\subset F_1\subset\ldots$ and $F_k\cdot F_l\subset F_{k+l}$
One checks that 
\begin{equation}
\label{eq-comm}
[F_k,F_l]\subset F_{k+l-1}
\end{equation}
so that associated graded algebra ${\rm Gr}(U(L,R))$ is commutative.
In fact, the Poincar\'e--Birkhoff--Witt isomorphism states that when $L$ is projective as an $R$-module, the symmetrization map
\[
{\rm Sym}_R(L)\to {\rm Gr}(U(L,R)),\qquad X_1\cdots X_k\mapsto \frac{1}{k!}\sum_{\sigma\in S_k}X_{\sigma(1)}\cdots X_{\sigma(k)}
\] 
is an isomorphism of algebras. In fact, following \cite{lsx}, it is possible to construct a lifting $\mathsf{pbw}:{\rm Sym}_R(L)\to U(L,R)$ of this map 
by choosing an $L$-connection $\nabla^L$, c.f.\ \cref{rk:con},  on $L$. With that, the map is defined inductively by
\begin{align}
\label{eq-pbw}
\mathsf{pbw}_k(\vec{X}):=&\frac{1}{k}\sum_{i=1}^kX_i\mathsf{pbw}_{k-1}(\vec{X}_{(i)})-\frac{1}{k}\sum_{i=1}^k\mathsf{pbw}_{k-1}(\nabla^L_{X_i}\vec{X}_{(i)}),
\end{align}
where $\nabla^L_{X_i}$ acts on ${\rm Sym}_R(L)$ by the tensor product connection.
\subsection{Extension to the adjoint complex}\label{subsect:extadj}
The aim of this section is to extend the PBW map to a morphism 
\[
\widetilde{\mathsf{pbw}}:\left({\rm Sym}({\rm ad}(L)),-\delta\right)\longrightarrow \left(C^\bullet(R,\sfU(L,R)),b\right)
\]
of cochain complexes. Its definition uses the cup product $\cup$ in \cref{eq-cup} on Hochschild cochains, together with the inclusion ${\rm Der}(R)\subset C^1(R,\sfU)$. 
For this we now choose a ${\rm Der}(R)$-connection $\nabla$ on $L$ and consider the associated basic connection $\nabla^b$ described in \cref{sec-ad} which acts on $\bigwedge^p{\rm Der}(R)\otimes_R{\rm Sym}^q_R(L)$ by
\begin{align*}
\nabla^b_Y(D_1\wedge\ldots\wedge D_p\otimes X_1\cdots X_q)&=\sum_{i=1}^pD_1\wedge\ldots\wedge \nabla^{{\rm Der}(R)}_Y(D_i)\wedge\ldots\wedge D_p\otimes X_1\cdots X_q\\
&+\sum_{j=1}^qD_1\wedge\ldots\wedge D_p\otimes X_1\cdots \nabla^L_{X}(X_j)\cdots X_q
\end{align*}

The extended PBW map is defined by a recursive formula similar to \cref{eq-pbw}:
\begin{align}
\label{eq-pbwp}
\notag \mathsf{\widetilde{pbw}}^{p,q}(\vec{D}  \otimes \vec{X}):=&
\sum_{i=1}^p(-1)^{i+1}D_i\cup \mathsf{\widetilde{pbw}}^{p-1,q}(\vec{D}_{(i)} \otimes \vec{X})\\
&+\sum_{i=1}^qX_i\mathsf{\widetilde{pbw}}^{p,q-1}(\vec{D} \otimes \vec{X}_{(i)})
 -\sum_{i=1}^q\widetilde{\mathsf{pbw}}^{p,q-1}\left(\nabla^b_{X_i}\left(\vec{D} \otimes \vec{X}_{(i)}\right)\right)
%\\ \notag &-\sum_{i=1}^p(-1)^i\mathsf{PBW}^{p,q-1}\left(\nabla_{D_i}\left(D_1\wedge\ldots\wedge\hat{D}_i\wedge\ldots\wedge D_p \otimes X_1,\ldots,X_q\right)\right)
\end{align}
Remark that the right hand side of the equation above contains only terms of total degree $p+q-1$, so that this formula indeed recursively defines a map $\widetilde{\mathsf{pbw}}$ by induction on the total degree. 
\begin{lemma}\label{lem-pbwcomplxmap}
The map $\widetilde{\mathsf{pbw}}$ is a morphism of complexes, that is:
\[\pbww^{p,q}\circ\delta =-b\circ\pbww^{p-1,q+1}.\]
\end{lemma}
\begin{proof}
This is done by induction on the total degree $p+q$. First remark that for $q=0$ we have 
\[
\widetilde{\mathsf{pbw}}^{p,0}(D_1\wedge\ldots\wedge D_p)(r_1,\ldots,r_p)=\sum_{\sigma\in S_p} D_{\sigma(1)}(r_1)\cup\cdots \cup D_{\sigma(p)}(r_p).
\]
This is the usual Hochschild--Kostant--Rosenberg map combined with the inclusion $R\subset\sfU(L,R)$, and is compatible with differentials.
Next, assume that 
\begin{align*}
\pbww^{p,q}\circ\delta&=-b\circ\pbww^{p-1,q+1}&\pbww^{p+1,q-1}\circ\delta&=-b\circ\pbww^{p,q}
\end{align*}
holds true for given $p$ and $q$.
With this we then compute
\begin{align*}
&\pbww^{p+1,q}(\delta(\vec{D}\otimes\vec{X}))=\sum_{i=1}^{q+1}\pbww^{p+1,q}(\rho(X_i)\wedge\vec{D}\otimes\vec{X_{(i)}})\\
&=\sum_{i\neq j} X_j\pbww^{p+1,q-1}(\rho(X_i)\wedge\vec{D}\otimes\vec{X_{(i,i)}})-\sum_{i\neq j}\pbww^{p+1,q-1}(\nabla^b_{X_j}(\rho(X_j)\wedge\vec{D}\otimes\vec{X_{(i,j)}}))\\
&+\sum_{i,j}(-1)^i D_i\cup\pbww^{p,q}(\rho(X_j)\wedge\vec{D_{(i)}}\otimes\vec{X_{(j)}})+\sum_{i=1}^{q+1}\rho(X_i)\cup\pbww^{p,q}(\vec{D}\otimes\vec{X_{(i)}})\\
&=\sum_{j=1}^{q+1}X_j\pbww^{p+1,q-1}(\delta(\vec{D}\otimes\vec{X_{(j)}}))-\sum_{j=1}^{q+1}\pbww^{p+1,q-1}(\nabla^b_{X_j}(\delta(\vec{D}\otimes\vec{X_{(j)}})))\\
&+\sum_{i=1}^{p}(-1)^i D_i\cup \pbww^{p,q}(\delta(\vec{D_{(i)}}\otimes\vec{X}))+\sum_{i=1}^{q+1}\rho(X_i)\cup\pbww^{p,q}(\vec{D}\otimes\vec{X_{(i)}})\\
&=-\sum_{j=1}^{q+1}X_j b(\pbww^{p,q}(\vec{D}\otimes\vec{X_{(j)}}))+\sum_{j=1}^{q+1}b(\pbww^{p,q}(\nabla^b_{X_j}(\vec{D}\otimes\vec{X_{(j)}}))\\
&-\sum_{i=1}^{p}(-1)^i D_i\cup b(\pbww^{p-1,q+1}(\vec{D_{(i)}}\otimes\vec{X}))-\sum_{i=1}^{q+1}\rho(X_i)\cup\pbww^{p,q}(\vec{D}\otimes\vec{X_{(i)}})\\
&=-\sum_{j=1}^{q+1}b(X_j\pbww^{p,q}(\vec{D}\otimes\vec{X_{(j)}}))+\sum_{j=1}^{q+1}b(\pbww^{p,q}(\nabla^b_{X_j}(\vec{D}\otimes\vec{X_{(j)}}))\\
&-\sum_{i=1}^p(-1)^{i+1}b(D_i\cup \pbww^{p-1,q}(\vec{D_{(i)}}\otimes\vec{X}))\\
&=-b(\pbww^{p,q+1}(\vec{D}\otimes\vec{X}))
\end{align*}
as required. Here we have used the identities
\[
[\delta,\nabla_X^b]=0,\qquad X(b\phi)=b(X\phi)-\rho(X)\cup\phi, \qquad b (D\cup \phi)=-D\cup b\phi.
\]

This proves the Lemma.
\end{proof}

\subsection{Extension to the nonlinear Chevalley--Eilenberg complex}
Before we proceed to extend the PBW map to the nonlinear Chevalley--Eilenberg complex, we introduce the following objects:
\begin{definition}
Let $\nabla$ be a connection on $L$ with associated basic connection $\nabla^b$. For $X,Y,Z\in L$ and $D\in{\rm Der}(R)$, define
\begin{align*}
\eta_Y(D,X)&:=[Y,\nabla_D(X)]-\nabla_{[\rho(Y),D]}(X)-\nabla_D[Y,X]\in L,\\
\eta^b_Y(X,D)&:=[Y,\nabla_X(D)]-\nabla_{[Y,X]}(D)-\nabla_X[\rho(Y),D]\in{\rm Der}(R),\\
\eta^b_Y(X,Z)&:=[Y,\nabla_X(Z)]-\nabla_{[Y,X]}(Z)-\nabla_X([Y,Z])\in L.
\end{align*}
\end{definition}
We have the relations
\[
\rho(\eta^b_Y(X,Z))=\eta^b_Y(X,\rho(Z)),\qquad \eta_Y(\rho(Z),X)=\eta^b_Y(X,Z),\qquad \rho(\eta_Y(D,X))=\eta^b_Y(X,D).
\]

The following properties follow from a straightforward computation:
\begin{subequations}
\begin{align}
\eta(rD,X)(Y)&=r\eta(D,X)(Y)\\
\eta(D,rX)(Y)&=r\eta(D,X)(Y)
\\
\eta(D,X)(rY)&=r\eta(D,X)(Y)\\
\notag & -\nabla_D(X)(r)Y+D(X(r))Y+X(f)\nabla_D(Y)-D(r)[Y,X].
\end{align}
\end{subequations}
We have the identity
\begin{align}
\label{eq-us}
\notag(\mathcal{L}_Y\circ \nabla_Z-\nabla_Z\circ \mathcal{L}_Y)(\vec{D}\otimes\vec{X})=&\nabla_{[Y,Z]}(\vec{D}\otimes\vec{X})\\
&+\sum_i(-1)^{i+1}\eta^b_Y(D_i,Z)\wedge\vec{D}_{(i)}\otimes\vec{X}+\sum_j\vec{D}\otimes\eta^b_Y(X_j,Z)\vec{X}_{(j)}.
\end{align}

With this form, we define maps
\[
F:\bigwedge^p{\rm Der}(R)\otimes_R{\rm Sym}^q_R(L)\to\bigwedge^{p-1}{\rm Der}(R)\otimes_R{\rm Sym}^{q+1}_R(L),
\]
given by
\[
F_Y(\vec{D}\otimes\vec{X}):=\sum_{i,j}(-1)^{i+1}\vec{D}_{(i)}\otimes\eta_Y(D_i,X_j)\vec{X}_{(j)}.
\]
From direct computations we infer the identities:
\begin{align}
\label{eq-usa}
(F_Y\circ\delta+\delta\circ F_Y)(\vec{D}\otimes\vec{X})&=\sum_{i,j}\vec{D}\otimes\eta_Y^b(X_i,X_j)\vec{X}_{(i,j)}+\sum_{i,j}(-1)^{i+1}\eta^b_Y(X_j,D_i)\wedge\vec{D}_{(i)}\otimes \vec{X}_{(j)}.
\end{align}
and
\begin{equation}
\label{eq-Fce}
[\mathcal{L}_{Y_1},F_{Y_2}]-[\mathcal{L}_{Y_2},F_{Y_1}]=F_{[Y_1,Y_2]}
\end{equation}
Now for any $n\geq 0$ and $\vec{Y}\in\Lambda^n L$ we introduce maps
\[s^{n,p,q}_{\vec{Y}}:\bigwedge^p_R{\rm Der}(R)\otimes_R{\rm Sym}^q_R(L)\to C^{p-n}(R,\sfU(L,R))\]
recursively by
\begin{align}
\label{eq-s}
s^{n,p,q}_{\vec{Y}}(\vec{D}\otimes\vec{X}):=&\sum_{i=1}^p (-1)^{i+1+n}D_i\cup s^{n,p-1,q}_{\vec{Y}}(\vec{D_{(i)}}\otimes\vec{X})\\
\notag&+\sum_{j=1}^qX_js^{n,p,q-1}_{\vec{Y}}(\vec{D}\otimes\vec{X_{(j)}})\\
\notag&-\sum_{j=1}^qs^{n,p,q-1}_{\vec{Y}}(\nabla^b_{X_j}(\vec{D}\otimes\vec{X_{(j)}}))\\
\notag&+\sum_{i=1}^n(-1)^{i+n}s^{n-1,p-1,q}_{\vec{Y}_{(i)}}(F_{Y_i}(\vec{D}\otimes\vec{X}))
\end{align}
Observe that $s^0=\pbww$ and that $s^{n,p,q}_{\vec{Y}}$ maps into $\Hom(R^{\otimes (p-n)},\ulr^{\leq q})$.

Clearly, $s^n_{\vec{Y}}$ is a morphism of $R$-modules, i.e., we have
\[
s^n_{\vec{Y}}\circ r=r\circ s^n_{\vec{Y}},\qquad\mbox{for all}~r\in R.
\]
\begin{proposition}
\label{prop-h}
For any $n\geq 0$ and $\vec{Y}\in L$, the following identity holds true:
\begin{equation}
\sum_{i=1}^{n+1}(-1)^{i+1} [\mathcal{L}_{Y_i},s^{n,p,q}_{\vec{Y}_{(i)}}]+\sum_{i<j}(-1)^{i+j} s^{n,p,q}_{[Y_i,Y_j],\vec{Y}_{(i,j)}}=b\circ s^{n+1,p,q}_{\vec{Y}}+(-1)^{n+1} s^{n+1,p+1,q-1}_{\vec{Y}}\circ\delta.
\end{equation}
\end{proposition}
\begin{proof}
We start with $n=q=0$. This case follows directly since
\[\pbww^{p,0}(D_1\wedge\cdots\wedge D_p)=\sum_{\sigma\in S_p}D_{\sigma(1)}\cup\cdots\cup D_{\sigma(p)}\]
which satisfies
\[[\mathcal{L}_Y,\pbww^{p,0}]=0.\]
On the other hand for $q=0$ we have $\delta=0$ and $s^{p,0}_Y=0$, so that the result holds.

We now do induction on $n+p+q$ by assuming
\begin{align*}
\sum_{i=1}^{n+1}(-1)^{i+1} [\mathcal{L}_{Y_i},s^{n,p-1,q}_{\vec{Y_{(i)}}}]+\sum_{i<j}(-1)^{i+j} s^{n,p-1,q}_{[Y_i,Y_j],\vec{Y_{(i,j)}}}&=b\circ s^{n+1,p-1,q}_{\vec{Y}}+(-1)^{n+1} s^{n+1,p,q-1}_{\vec{Y}}\circ\delta,\\
\sum_{i=1}^{n+1}(-1)^{i+1} [\mathcal{L}_{Y_i},s^{n,p,q-1}_{\vec{Y_{(i)}}}]+\sum_{i<j}(-1)^{i+j} s^{n,p,q-1}_{[Y_i,Y_j],\vec{Y_{(i,j)}}}&=b\circ s^{n+1,p,q-1}_{\vec{Y}}+(-1)^{n+1} s^{n+1,p+1,q-2}_{\vec{Y}}\circ\delta,
\end{align*}
and
\begin{equation*}
\sum_{i=1}^{n}(-1)^{i+1} [\mathcal{L}_{Y_i},s^{n-1,p,q-1}_{\vec{Y_{(i)}}}]+\sum_{i<j}(-1)^{i+j} s^{n-1,p,q-1}_{[Y_i,Y_j],\vec{Y_{(i,j)}}}=b\circ s^{n,p,q-1}_{\vec{Y}}+(-1)^{n} s^{n,p+1,q-2}_{\vec{Y}}\circ\delta.
\end{equation*}

In the case that $n=0$ we read the last assumption as
\[b\circ s^{0,p,q-1}+s^{0,p+1,q-2}\circ\delta=0,\]
which is true by \cref{lem-pbwcomplxmap} and the observation that $s^0=\pbww$.

Using the recursive definition of $s^{n,p,q}$ we find after some manipulations:
\begin{align*}
\sum_{i=1}^{n+1}(-1)^{i+1}[\mathcal{L}_{Y_i},s^{n,p,q}_{\vec{Y_{(i)}}}](\vec{D}\otimes\vec{X})=&\sum_{i=1}^{n+1}\sum_{j=1}^p(-1)^{i+j+n}D_j\cup [\mathcal{L}_{Y_i},s^{n,p-1,q}_{\vec{Y_{(i)}}}](\vec{D_{(j)}}\otimes\vec{X})\tag{a.1}\\
&+\sum_{i=1}^{n+1}\sum_{j=1}^q(-1)^{i+1}X_j[\mathcal{L}_{Y_i},s^{n,p,q-1}_{\vec{Y_{(i)}}}](\vec{D}\otimes\vec{X_{(j)}}))\tag{a.2}\\
&+\sum_{i=1}^{n+1}\sum_{j=1}^q(-1)^i[\mathcal{L}_{Y_i},s^{n,p,q-1}_{\vec{Y_{(i)}}}](\nabla^b_{X_j}(\vec{D}\otimes\vec{X_{(j)}})\tag{a.3}\\
&+\sum_{i=1}^{n+1}\sum_{j=1}^q(-1)^is^{n,p,q-1}_{\vec{Y_{(i)}}}(\mathcal{L}_{Y_i}(\nabla^b_{X_j}(\vec{D}\otimes\vec{X_{(j)}})))]\tag{a.4.1}\\
&-\sum_{i=1}^{n+1}\sum_{j=1}^q(-1)^i s^{n,p,q-1}_{\vec{Y_{(i)}}}(\nabla^b_{[Y_i,X_j]}(\vec{D}\otimes\vec{X_{(j)}}))\tag{a.4.2}\\
&-\sum_{i=1}^{n+1}\sum_{j=1}^q(-1)^i s^{n,p,q-1}_{\vec{Y_{(i)}}}(\nabla^b_{[Y_i,X_j]}(\mathcal{L}_{Y_i}(\vec{D}\otimes\vec{X_{(j)}}))\tag{a.4.3}\\
&+\sum_{j<i}^{n+1}(-1)^{i+j+n+1}[\mathcal{L}_{Y_i},s^{n-1,p-1,q}_{\vec{Y_{(j,i)}}}](F_{Y_j}(\vec{D}\otimes\vec{X}))\tag{a.5}\\
&+\sum_{j<i}^{n+1}(-1)^{i+j+n+1}s^{n-1,p-1,q}_{\vec{Y_{(j,i)}}}([\mathcal{L}_{Y_i},F_{Y_j}](\vec{D}\otimes\vec{X}))\tag{a.6}\\
&+\sum_{i<j}^{n+1}(-1)^{i+j+n}[\mathcal{L}_{Y_i},s^{n-1,p-1,q}_{\vec{Y_{(j,i)}}}](F_{Y_j}(\vec{D}\otimes\vec{X}))\tag{a.7}\\
&+\sum_{i<j}^{n+1}(-1)^{i+j+n}s^{n-1,p-1,q}_{\vec{Y_{(j,i)}}}([\mathcal{L}_{Y_i},F_{Y_j}](\vec{D}\otimes\vec{X}))\tag{a.8}
\end{align*}
We note that by \eqref{eq-us} and \eqref{eq-usa} we have
\[
\text{(a.4.1)}+\text{(a.4.2)}+\text{(a.4.3)}=\sum_{i=1}^{n+1}(-1)^is^{n,p,q}_{\vec{Y_{(i)}}}(\{F_{Y_i},\delta\}(\vec{D}\otimes\vec{X}))=\text{(a.4)}\]

Next, we have
\begin{align*}
\sum_{i<j}^{n+1}(-1)^{i+j}s^{n,p,q}_{[Y_i,Y_j],\vec{Y_{(i,j)}}}(\vec{D}\otimes\vec{X})=&\sum_{i<j}^{n+1}\sum_{k=1}^p(-1)^{i+j+k+1+n}D_k\cup s^{n,p-1,q}_{[Y_i,Y_j],\vec{Y_{(i,j)}}}(\vec{D_{(k)}}\otimes\vec{X})\tag{b.1}\\
&+\sum_{i<j}^{n+1}\sum_{k=1}^q (-1)^{i+j}X_k s^{n,p,q-1}_{[Y_i,Y_j],\vec{Y_{(i,j)}}}(\vec{D}\otimes\vec{X_{(k)}})\tag{b.2}\\
&+\sum_{i<j}^{n+1}\sum_{k=1}^q(-1)^{i+j+1} s^{n,p,q-1}_{[Y_i,Y_j],\vec{Y_{(i,j)}}}(\nabla^b_{X_k}(\vec{D}\otimes\vec{D_{(k)}}))\tag{b.3}\\
&+\sum_{i<j}^{n+1}(-1)^{i+j+n+1}s^{n-1,p,q-1}_{\vec{Y_{(i,j)}}}(F_{[Y_i,Y_j]}(\vec{D}\otimes\vec{X}))\tag{b.4}\\
&+\sum_{k<i<j}^{n+1}(-1)^{i+j+k+n+1}s^{n-1,p,q-1}_{[Y_i,Y_j],\vec{Y_{(k,i,j)}}}(F_{Y_k}(\vec{D}\otimes\vec{X}))\tag{b.5}\\
&+\sum_{i<k<j}^{n+1}(-1)^{i+j+k+n}s^{n-1,p-1,q}_{[Y_i,Y_j],\vec{Y_{(i,k,j)}}}(F_{Y_k}(\vec{D}\otimes\vec{X}))\tag{b.6}\\
&+\sum_{i<j<k}^{n+1}(-1)^{i+j+k+n+1}s^{n-1,p-1,q}_{[Y_i,Y_j],\vec{Y_{(i,j,k)}}}(F_{Y_k}(\vec{D}\otimes\vec{X}))\tag{b.7}
\end{align*}
By \eqref{eq-Fce}, the terms (a.6), (a.8) and (b.4) cancel against each other. Now we use the induction hypothesis to obtain
\begin{align*}
\text{(a.1)}+\text{(b.1)}
%=&\sum_{j=1}^p(-1)^{j+n+1}D_j\cup b(s^{n+1,p-1,q}_{\vec{Y}}(\vec{D_{(j)}}\otimes\vec{X}))\\
%&+\sum_{j=1}^p(-1)^j D_j\cup s^{n+1,p,q-1}_{\vec{Y}}(\delta(\vec{D_{(j)}}\otimes\vec{X}))\\
=&\sum_{j=1}^p(-1)^{j+n} b(D_j\cup s^{n+1,p-1,q}_{\vec{Y}}(\vec{D_{(j)}}\otimes\vec{X}))\\
&+\sum_{j=1}^p(-1)^j D_j\cup s^{n+1,p,q-1}_{\vec{Y}}(\delta(\vec{D_{(j)}}\otimes\vec{X})),
\end{align*}
as well as 
\begin{align*}
\text{(a.2)}+\text{(b.2)}
%=&\sum_{j=1}^q X_j b(s^{n+1,p,q-1}_{\vec{Y}}(\vec{D}\otimes\vec{X_{(j)}}))\\
%&+\sum_{j=1}^q (-1)^{n+1}X_j s^{n+1,p+1,q-2}_{\vec{Y}}(\delta(\vec{D}\otimes\vec{X_{(j)}}))\\
=&\sum_{j=1}^q b(X_js^{n+1,p,q-1}_{\vec{Y}}(\vec{D}\otimes\vec{X_{(j)}}))
+\sum_{j=1}^q \rho(X_j)\cup s^{n+1,p,q-1}_{\vec{Y}}(\vec{D}\otimes\vec{X_{(j)}})\\
&+\sum_{j=1}^q (-1)^{n+1}X_j s^{n+1,p+1,q-2}_{\vec{Y}}(\delta(\vec{D}\otimes\vec{X_{(j)}})),
\end{align*}
and
\begin{align*}
\text{(a.3)}+\text{(b.3)}
%=&-\sum_{j=1}^q b(s^{n+1,p,q-1}_{\vec{Y}}(\nabla^b_{X_j}(\vec{D}\otimes\vec{X_{(j)}})))\\
%&-\sum_{j=1}^q(-1)^{n+1} s^{n+1,p+1,q-2}_{\vec{Y}}(\delta(\nabla^b_{X_j}(\vec{D}\otimes\vec{X_{(j)}}))\\
=&-\sum_{j=1}^q b(s^{n+1,p,q-1}_{\vec{Y}}(\nabla^b_{X_j}(\vec{D}\otimes\vec{X_{(j)}})))\\
&-\sum_{j=1}^q(-1)^{n+1} s^{n+1,p+1,q-2}_{\vec{Y}}(\nabla^b_{X_j}(\delta(\vec{D}\otimes\vec{X_{(j)}}))),
\end{align*}
and finally
\begin{align*}
\text{(a.5)}+\text{(a.7)}+\text{(b.5)}+\text{(b.6)}+\text{(b.7)}=&\sum_{i=1}^{n+1}(-1)^{i+n+1} b(s^{n,p-1,q}_{\vec{Y_{(i)}}}(F_{Y_i}(\vec{D}\otimes\vec{X})))\\
&+\sum_{i=1}^{n+1}(-1)^{i+1}s^{n,p,q-1}_{\vec{Y_{(i)}}}(\delta(F_{Y_i}(\vec{D}\otimes\vec{X})))
\end{align*}
Combining all we obtain the desired equality stated in the proposition. 
\end{proof}

Now the maps $s$ induces maps:
\[ 
\mathfrak{s}^k: {\rm Hom}\left(\bigwedge^\bullet L,\bigwedge^\bullet_R{\rm Der}(R)\otimes{\rm Sym}_R(L)\right)\to{\rm Hom}\left(\bigwedge^{\bullet+k}L,C^{\bullet-k}_{\rm Hoch}(R,\sfU(L,R))\right)
\]
which takes an element $\phi\in {\rm Hom}\left(\bigwedge^i L,\bigwedge^p_R{\rm Der}(R)\otimes{\rm Sym}_R(L)\right)$ and maps it to the element $\mathfrak{s}^k\phi\in\Hom(\Lambda^{i+k}L,C^{p-k}_{\rm Hoch}(R,\sfU(L,R)))$ given by
\begin{align*}
(\mathfrak{s}^k\phi)(Y_1,...,Y_{k+i})=\frac{1}{k!i!}\sum_{\sigma\in S_{k+i}} s^k_{Y_{\sigma(1)},...,Y_{\sigma(k)}}(c(Y_{\sigma(k+1)},...,Y_{\sigma(k+i)})).
\end{align*}
\cref{prop-h} then immediately implies:
\begin{theorem}
\label{thm-epbw}
The map
\[\Phi: (\phi_0,...,\phi_n)\mapsto (\sum_{i=0}^n \mathfrak{s}^i(\phi_i),\sum_{i=0}^n\mathfrak{s}^i(\phi_{i+1}),...,\mathfrak{s}^0(\phi_{n-1})+\mathfrak{s}^1(\phi_n),\mathfrak{s}^0(\phi_n))\]
defines a morphism of cochain complexes:
\[\Phi: \sfOmega_{{\rm nl}}({\rm Sym}_R({\rm ad}(L,R)),\delta_{\text{CE}}\pm -\delta)\to\sfOmega(C^\bullet(R,\ulr),\delta_{\text{CE}}\pm b).\]
\end{theorem}

\begin{remark}\label{rmk-thm-epbw}
Considering low degrees, one easily observes that in general
\[
\widetilde{\mathsf{pbw}}\circ h_{r,Y}-h_{r,Y}\circ\widetilde{\mathsf{pbw}}\not =r\circ s^1_Y+s^1_{rY}.
\]
However, equalities of this kind would mean that the extended Poincar\'e--Birkhoff--Witt morphism maps the nonlinear Chevalley--Eilenberg cochains of ${\rm Sym}({\rm ad}(L))$ to those of $C^\bullet(R,\sfU)$. This is the reason why we can't write $\sfOmega_{\rm nl}(C^\bullet(R,\ulr)$ on the right hand side of \cref{thm-epbw}.
\end{remark}
\section{Hochschild cohomology}
\label{sec-hc}
We now come to the main result of the paper: the computation of the Hochschild cohomology of the universal enveloping algebra $\sfU(L,R)$ of a Lie--Rinehart algebra.
Our strategy is to first reduce the computation to the nonlinear Chevalley--Eilenberg complex of the quasi-module structure on the Hochschild complex $C^\bullet(R,\sfM)$.
We then use the extended PBW-morphism of the previous section to find the result computed in terms of the adjoint representation.
\subsection{A spectral sequence}
Let $(L,R)$ be a Lie--Rinehart pair and $\sfM$ a $\sfU(L,R)$-bimodule.
Here we recall the construction of a spectral sequence in \cite{KL} converging to the Hochschild cohomology $HH^\bullet(\sfU(L,R),\sfM)$.
By definition the Hochschild cohomology groups $HH^\bullet(\sfU(L,R),\sfM)$ are given as the right derived functors of
\[
{\rm Hom}_{\sfU^e}(\sfU,-): {}_{\sfU^e}\mathsf{Mod}\longrightarrow \mathsf{Vect}_{\mathbb{K}}.
\]
Consider now, as in \cite{llm} the functor 
\[
G:{}_{\sfU^e}\mathsf{Mod}\longrightarrow {}_{\sfU}\mathsf{Mod},\qquad G(\sfM):={\rm Hom}_{R^e}(R,M),
\]
where $(L,R)$ acts on $G(M)$ by
\begin{subequations}
\begin{align}
\label{RI1}
r_1\cdot\phi(r_2)&:=r_1\phi(r_2),\\
\label{RI2}
(X\cdot\phi)(r)&:=[X,\phi(r)]-\phi(X(r)).
\end{align}
\end{subequations}
It is proved in \cite[Prop. 3.4.1.]{llm} that this functor is right adjoint to the induction functor
\begin{equation}
\label{eq-ind}
I:\mathsf{Mod}_\sfU\to\mathsf{Mod}_{\sfU^e},\qquad \sfM\mapsto \sfU\otimes_R\sfM,
\end{equation}
 where the left $\sfU$-module structure is the obvious one and $\sfU$ acts from the right as
\[
(\sfD\otimes m)\cdot X:=\sfD X\otimes m-\sfD\otimes Xm,\qquad \mbox{for}~ \sfD\in\sfU,~m\in \sfM~\mbox{and}~X\in L. 
\] 
We denote by $F: \mathsf{Mod}_{\sfU}\to \mathsf{Vect}_{\mathbb{K}}$ the functor $F(N):={\rm Hom}_{\sfU}(R,N)$ taking $\sfU(L,R)$-invariants.
We then have
\begin{lemma}
There is a natural isomorphism of functors
\[
{\rm Hom}_{\sfU^e}(\sfU,-)\cong F\circ G.
\]
\end{lemma} 
\begin{proof}
Let $\sfM\in {}_{\sfU^e}\mathsf{Mod}$ be a $\sfU(L,R)$-bimodule. Using the adjunction $I\dashv G$ mentioned above we see that
\begin{align*}
F(G(\sfM))&={\rm Hom}_\sfU(R,G(\sfM))\cong {\rm Hom}_{\sfU^e}(I(R),\sfM)={\rm Hom}_{\sfU^e}(\sfU,\sfM). 
\end{align*}
This isomorphism is natural in $\sfM$, and therefore proves the Lemma.
%It is easy to check that $F(G(M))=[M,U]$ for the $(L,R)$-action given by \cref{RI1,RI2} on $G(M)$.
\end{proof}
Because of this Lemma there is a Grothendieck spectral sequence
\[
E_2^{p,q}=H^p_{R}(L,H^q(R,\sfM))\quad\Longrightarrow\quad HH^{p+q}(\sfU(L,R),\sfM).
\]
This is exactly the spectral sequence constructed in \cite[Cor. 3.3]{KL}.
Explicitly this spectral sequence can be constructed by means of a double complex by choosing an injective resolution $I^\bullet\to\sfM$ of $\sfM$ 
in the category of $\sfU^e$-modules and a projective resolution $Q_\bullet\to R$ of $R$ in the category of $\sfU$-modules. With this we set
\begin{equation}
\label{eq-dch}
\mathsf{X}^{\bullet,\bullet}:={\rm Hom}_\sfU(Q_\bullet,G(I^\bullet)),
\end{equation}
with the differentials induced those from $I^\bullet$ and $Q_\bullet$.

%\begin{proposition}
%There exists a spectral sequence with
%\[
%E_2^{p,q}\cong H^{p+q}\left(\bigwedge^\bullet_R{\rm Der}(R)\otimes {\rm Sym}_R(L),\delta_{\rm Koszul}\right) \quad \Longrightarrow \quad H^q(R,U(L,R))
%\]
%\end{proposition}
%\begin{proof}
%By the Poincar\'e--Birkhoff--Witt (PBW) isomorphism we have $U(L,R)\cong{\rm Sym}_R(L)$, and we have
%\[
%E_1^{p,q}\cong H^{p+q}(R,{\rm Sym}_R(L))\cong\bigwedge^{p+q}_R{\rm Der}(R)\otimes_R{\rm Sym}_R(L).
%\]
%\end{proof}
\subsection{Reduction to the nonlinear Chevalley--Eilenberg complex}
In this section we reduce the computation of the Hochschild cohomology $H^\bullet(\sfU,\sfM)$ to the cohomology of the nonlinear Chevalley--Eilenberg complex
$\sfOmega_{\rm nl}(L;\sfM)$. We start with an example in degree $1$:
\begin{example}
In degree $1$, a cocycle in $C^1(\sfU,\sfM)$ is the same as a derivation $\Delta:\sfU\to\sfM$.
On the other hand, a cocycle of degree $1$ in $\sfOmega_{\rm nl}(L;C^\bullet(R,\sfM))$ is given by a pair  $\phi=(\phi_0,\phi_1)$
with $\phi_0\in {\rm Hom}(R,\sfM)$ and $\phi_{1}\in {\rm Hom}(L,\sfM)$ satisfying
\[
\phi_{1}(rX)=r\phi_{1}(X)+\phi_0(r)X.
\]
Imposing the cocycle condition leads to the equations
\begin{align*}
[X,\phi_1(Y)]-[Y,\phi_1(X)]+\phi_{1}([X,Y])]&=0,\\
[X,\phi_{0}(r)]+\phi_0(X(r))&=[r,\phi_{1}(X)],\\
r_1\phi_{0}(r_2)-\phi_{0}(r_1r_2)+\phi_{0}(r_1)r_2&=0.
\end{align*}
Now, given a derivation $\Delta:\sfU\to\sfM$, it is easy to check that 
\[
\phi_0(r):=\Delta(r),\qquad\phi_1(X):=\Delta(X),
\]
satisfy the four equations above. Conversely, given $(\phi_0,\phi_1)$ satisfying the equations above, it is straightforward to check that extending by imposing the derivation property
\[
\Delta(rX_1\cdots X_k):=\phi_0(r)X_1\cdots X_k+\sum_{i=1}^krX_1\cdots X_{i-1}\Delta(X_i) X_i\cdots X_k,
\]
the equations above ensure that $\Delta$ is well-defined. 
\end{example}
In the higher degrees the relation is not so straightforward. Nevertheless we have:
\begin{proposition}
\label{prop-nlh}
Let $(L,R)$ be a Lie--Rinehart algebra and $\sfM$ a $\sfU(L,R)$-bimodule. 
The cohomology of the complex
\[
\left(\sfOmega_{\rm nl}(L;C^\bullet(R,\sfM)),\delta_{\mathsf{CE}}\pm b\right)
\] 
equals the Hochschild cohomology $HH^\bullet(\sfU(L,R),\sfM)$.
\end{proposition}
\begin{proof}
Let $I^\bullet\to\sfM$ be an injective resolution of $\sfU^e$-bimodules. We have that 
\[
{\rm Hom}_{R^e}(R,I^\bullet)\longrightarrow {\rm Tot}({\rm Hom}_{R^e}(\mathsf{Bar}_\bullet(R),I^\bullet))\longleftarrow {\rm Hom}_{R^e}(\mathsf{Bar}_\bullet(R),\sfM)
\]
are quasi-isomorphisms, as we can see by filtering the middle complex by either its rows 
or its columns. The complex on the right hand side is simply ${\rm Hom}_{R^e}(\mathsf{Bar}_\bullet(R),\sfM)\cong C^\bullet(R,\sfM)$ with its quasi module structure given by \cref{ar,acl,eq-hqa}.
The complex on the left hand side carries a canonical $(L,R)$-module structure given by \cref{RI1,RI2}. We now combine these two $R$ and $L$ module structures to define a 
quasi $(L,R)$-module structure on the middle complex ${\rm Tot}({\rm Hom}_{R^e}(\mathsf{Bar}_\bullet(R),I^\bullet))$ turning both maps into morphisms of quasi $(L,R)$-modules as in \cref{prop-invqi}. By that proposition we therefore see that the two complexes
\[
\sfOmega_{\rm nl}(L;C^\bullet(R,\sfM))\qquad\sfOmega(L;G(I^\bullet))
\]
are quasi-isomorphic. But the right hand complex is exactly the total complex in \cref{eq-dch} 
for the Chevalley--Eilenberg resolution $\sfU(L,R)\otimes_R\bigwedge^\bullet_RL\to R$ of $R$ in the category of $\sfU(L,R)$-modules. It is shown in \cite[Thm 3.2.]{KL} that the cohomology of
this complex is isomorphic to the Hochschild cohomology $HH^\bullet(\sfU(L,R),\sfM)$. 
\end{proof}

\subsection{Relation with the adjoint representation}
We now take $\sfM=\sfU(L,R)$ and consider the complex 
$\sfOmega^\bullet_{\rm nl}(L;C^\bullet(R,\sfU(L,R)))$. By \cref{prop-nlh}, this complex 
computes the Hochschild cohomology $HH^\bullet(\sfU(L,R),\sfU(L,R))$.
The main result is now:
\begin{proposition}
\label{prop:pc}
The cohomology of the complex $\left(\sfOmega^\bullet_{\rm nl}(L;C^\bullet(R,\sfU(L,R))),\delta_{\mathsf{CE}}\pm b\right)$ is isomorphic to the Poisson cohomology $H_{\rm Pois}({\rm Sym}_R(L))$.
\end{proposition}
\begin{remark}
Combined with \cref{prop-nlh}, this completes the proof of \cref{thm-m1}.
\end{remark}
\begin{proof}
In \cref{thm-epbw} we have constructed a morphism
\[
\Phi:\sfOmega(L;{\rm Sym}_R({\rm ad}(L,R))\longrightarrow \sfOmega(L;C^\bullet(R,\sfU)).
\]
As remarked in \cref{rmk-thm-epbw} this morphism does not preserve non-linear chains. 
We therefore define
\[
\mathcal{H}(L,R):= \Phi(\sfOmega_{\rm nl}(L;{\rm Sym}_R({\rm ad}(L,R)))+\sfOmega_{\rm nl}(L;C^\bullet(R,\sfU))\subset \sfOmega(L;C^\bullet(R,\sfU)).
\]
Below we shall  argue that the maps
\begin{equation}\label{eq-zigzag}
\begin{tikzcd}
&\mathcal{H}(L,R)&\\
\sfOmega_{\rm nl}(L;C^\bullet(R,\sfU))\arrow[ru, hook] &&\sfOmega_{\rm nl}(L;{\rm Sym}_R({\rm ad}(L,R)))\arrow[lu, "\Phi"', hook']
\end{tikzcd}
\end{equation}
are quasi-isomorphisms.

To see this, we introduce filtrations on $\sfOmega(L;{\rm Sym}_R({\rm ad}(L,R))$ and $\sfOmega(L;C^\bullet(R,\sfU))$ and in turn also on the non-linear subcomplexes. The filtrations will use the filtration on $\sfU(L,R)$ (see \cref{section:pbw}) and the grading on ${\rm Sym}_RL$. The filtration on $\sfOmega(L;{\rm Sym}_R({\rm ad}(L,R))$ is given by
\begin{multline}F_k(\sfOmega^n(L;{\rm Sym}_R({\rm ad}(L,R)))=\{(\phi_0,...,\phi_n)\in \sfOmega^n(L;{\rm Sym}_R({\rm ad}(L,R)):\\ \phi_i(\Lambda^{n-i}L)\subset \Lambda_R^i\Der(R)\otimes{\rm Sym}^{\leq k+n-i}_RL\},
\end{multline}
while on $\sfOmega(L;C^\bullet(R,\sfU))$ it is given by
\begin{multline}F_k(\sfOmega^n(L;C^\bullet(R,\sfU)))=\{(\phi_0,...,\phi_n)\in\sfOmega^n(L;C^\bullet(R,\sfU)):\\
\phi_i(\Lambda^{n-i}L)(R^{\otimes i})\subset F_{k+n-i}(\sfU(L,R))\}.
\end{multline}
There are a few remarks to make here:
\begin{itemize}
\item The filtration on $\sfOmega(L;{\rm Sym}_R({\rm ad}(L,R))$ is induced by the obvious grading;
\item The differential on $\sfOmega(L;{\rm Sym}_R({\rm ad}(L,R))$ lowers the graded degree by $1$;
\item Using the Poincar\'e-Birkhoff-Witt map $\mathsf{pbw}$ we see that the graded quotient of $\sfOmega(L;C^\bullet(R,\sfU))$ is isomorphic to $\sfOmega(L;C^\bullet(R,{\rm Sym}_RL))$, and using the Hochschild-Kostant-Rosenberg Theorem we see that the cohomology of the graded quotient of $\sfOmega(L;C^\bullet(R,\sfU))$ is isomorphic to $\sfOmega(L;{\rm Sym}_R({\rm ad}(L,R))$.
\item Since the Poisson bracket on ${\rm Sym}_RL$ is induced by the filtration on $\sfU(L,R)$ and the PBW-morphism, the induced differential on the cohomology of the graded quotient of $\sfOmega(L;C^\bullet(R,\sfU))$ equals the differential on $\sfOmega(L;{\rm Sym}_R({\rm ad}(L,R))$.
\end{itemize}
In particular we see that if we write out the spectral sequences for the filtered complexes $\sfOmega(L;{\rm Sym}_R({\rm ad}(L,R))$ and $\sfOmega(L;C^\bullet(R,\sfU))$, they are isomorphic to each other on the second page.

Next, we investigate how the zigzag \cref{eq-zigzag} interacts with this filtration. Here we make two remarks:
\begin{itemize}
\item The map $\widetilde{\mathsf{pbw}}$ from \cref{subsect:extadj} satisfies
\[\widetilde{\mathsf{pbw}}(\overrightarrow{D}\otimes\overrightarrow{X})(r_1,...,r_p)=\sum_{\sigma\in S_p}D_{\sigma(1)}(r_1)\cdots D_{\sigma(p)}(r_p)\mathsf{pbw}(\overrightarrow{X})+{\rm L.O.T.}\]
\item The homotopies $h_{r,X}$ on ${\rm Sym}_R({\rm ad}(L,R))$ and $C^\bullet(R,\sfU)$ are compatible in the following way:
\[h_{r,X}\circ\widetilde{\mathsf{pbw}}=\widetilde{\mathsf{pbw}}\circ h_{r,X}+{\rm L.O.T.}\]
\end{itemize}
In particular we see that not only are the graded quotients of $\sfOmega(L;{\rm Sym}_R({\rm ad}(L,R))$ and $\sfOmega(L;C^\bullet(R,\sfU))$ quasi-isomorphic as chain-complexes, the same actually holds true for $\sfOmega_{\rm nl}(L;{\rm Sym}_R({\rm ad}(L,R))$ and $\sfOmega_{\rm nl}(L;C^\bullet(R,\sfU))$. On top of this, the maps:
\begin{align*}
\sfOmega_{\rm nl}(L;C^\bullet(R,\sfU))&\hookrightarrow \sfOmega(L;C^\bullet(R,\sfU)),&\sfOmega_{\rm nl}(L;{\rm Sym}_R({\rm ad}(L,R))&\xrightarrow{\Phi}\sfOmega(L;C^\bullet(R,\sfU))
\end{align*}
induce the following commutative triangle on the level of graded quotients:
\[
\begin{tikzcd}
 &{\rm Gr}(\sfOmega(L;C^\bullet(R,\mathsf{U})))&                                                                              \\
{\rm Gr}(\mathsf{\Omega}_{\rm nl}(L;C^\bullet(R,\sfU))) \arrow[ru, hook] \arrow[rr, "\simeq", dotted] &                                                  
& {\rm Gr}(\mathsf{\Omega}_{\rm nl}(L;{\rm Sym}_R({\rm ad}(L,R))) \arrow[lu, "\Phi"', hook']
\end{tikzcd}
\]
where the horizontal quasi-isomorphism is the one described above.
Combining all this, we see that the zigzag \cref{eq-zigzag} induces quasi-isomorphisms on the level of graded quotients, and hence the zigzag itself is composed of quasi-isomorphisms.
\end{proof}
\section{Hochschild and cyclic homology}
\label{sec-hm}
\subsection{Duality for Lie--Rinehart algebras}
Recall that an algebra $A$ is said to satisfy {\em van de Bergh duality} in dimension $n$ if
\begin{itemize}
\item[$i)$] $A$ is homologically smooth, i.e., it has a bounded resolution by finitely generated projective $A^e$-modules,
\item[$ii)$] we have that
\[
{\rm Ext}^\bullet_{A^e}(A,A^e)=\begin{cases} \omega_A,&\bullet=n,\\0,&\bullet \not =0,\end{cases}
\]
and the $A^e$-bimodule $\omega_A$ is invertible.
\end{itemize}
If this is the case there are functorial isomorphisms \cite{vdb}
\begin{equation}
\label{eq:vdblr}
HH_p(A,M)\cong HH^{n-p}(A,\omega_A^{-1}\otimes_A M).
\end{equation}

Let $(L,R)$ be a Lie--Rinehart algebra and assume now, in addition to $R$ being smooth and $L$ projective, that $R$ satisfies van den Bergh duality in dimension $n$.
Using the HKR-isomorphism this implies that 
\[
{\rm Ext}^n_{R^e}(R,R^e)\cong\bigwedge^n_R{\rm Der}(R)
\]
is invertible with inverse $\Omega^n_R:={\rm Hom}_R(\bigwedge^n_R{\rm Der}(R),R)$. If we furthermore assume that $L$ is finitely generated of constant rank $d$,
it is proved in \cite[Thm. 1]{llm} that $\sfU(L,R)$ satisfies van den Bergh duality in dimension $n+d$ with
\[
{\rm Ext}_{\sfU^e}^{n+d}(\sfU,\sfU^e)\cong I\left(\bigwedge_R^dL^\vee\otimes_R\bigwedge^n_R{\rm Der}(R)\right),
\]
where we recall the induction functor $I$ in \cref{eq-ind}. Here $\bigwedge_R^dL^\vee$ is the dualizing right $(L,R)$-module of \cite[Thm. 2.10]{huebschmann-dual},
and the tensor product 
\[
Q_L^\vee:=\bigwedge_R^dL^\vee\otimes_R\bigwedge^n_R{\rm Der}(R),
\]
is the (dual) $(L,R)$-module controlling the so-called modular class in  the 
geometric setting \cite{elw}. 
\subsection{Hochschild homology}
We now use \cref{eq:vdblr} to compute the Hochschild homology of $\sfU(L,R)$. By the PBW-isomorphism and \cref{lemma-ses} we have isomorphisms
\[
I(Q_L^\vee)\cong{\rm Sym}_R(L)\otimes_R\bigwedge_R^dL^\vee\otimes_R\bigwedge^n_R{\rm Der}(R)\cong\bigwedge^{n+d}_R\left({\rm Der}\left({\rm Sym}_R(L)\right)\right).
\]
We now claim that the entire \cref{sec:var-pbw} generalizes verbatim to produce a morphism
\begin{equation}
\label{eq-twm}
\Phi: \sfOmega_{\text{\rm nl}}({\rm Sym}_R(\text{ad}(L,R))\otimes_RQ^\vee_L,\delta_{\text{CE}}+\delta)\to\sfOmega(C^\bullet(R,I(Q^\vee_L)),\delta_{\text{CE}}+b).
\end{equation}
To see this, remark:
\begin{itemize}
\item The tensor product ${\rm Sym}_R(L)\otimes_RQ_L^\vee$ has a natural $L$ connection given by the tensor product connection given by $\nabla^L$ on ${\rm Sym}_R(L)$ and the $(L,R)$-module structure on $Q_L^\vee$. With this connection, the PBW morphism is defined by the same formula as in \cref{eq-pbw}.
\item We extend this connection to $\bigwedge^p{\rm Der}(R)\otimes_R{\rm Sym}^q_R(L)\otimes_RQ_L^\vee$ and define the extended PBW morphism by the 
same formula as in \cref{eq-pbwp} with the remark that the cup product with $D_i\in{\rm Der}(R)$ still makes sense on $C^\bullet(R,I(Q_L^\vee))$. \cref{lem-pbwcomplxmap}
still holds true in this context because the identities at the end of its proof continue to hold.
\item The homotopies $s$ are defined by the same \cref{eq-s} using the extended connection above and \cref{prop-h} continues to hold. 
\end{itemize}
The morphism \cref{eq-twm}, together with \cref{prop-nlh}, then gives the following twisted analogue of \cref{prop:pc}:
\begin{theorem}
Let $(L,R)$ be a Lie--Rinehart algebra with $L$ projective over $R$, and let $\sfM$ be an $(L,R)$-module.
There is a canonical isomorphism 
\[
H^\bullet(\sfU,I(\sfM))\cong H^\bullet({\rm Sym}({\rm ad})\otimes\sfM).
\]
\end{theorem}
Applying van den Bergh duality for Lie--Rinehart algebras we therefore conclude that
\begin{equation}
\label{eq-hi}
HH_\bullet(\sfU(L,R))\cong H^{n+d-\bullet}({\rm Sym}({\rm ad})\otimes Q_L)
\end{equation}

\subsection{Cyclic homology}
Let us take a closer look at the complex computing the cohomology in \cref{eq-hi}. 
In \cref{sec-ad} we have seen that the nonlinear Chevally--Eilenberg cochains associated to ${\rm Sym}({\rm ad})$ can be identified with the Poisson complex
$\mathfrak{X}^\bullet_{{\rm Sym}_R(L)}$. Now consider the short exact sequence
\[
0\longrightarrow {\rm Sym}_R(L)\otimes_R\Omega^1_R\longrightarrow \Omega^1_{{\rm Sym}_R(L)}\longrightarrow {\rm Sym}_R(L)\otimes_RL\longrightarrow 0,
\]
dual to the short exact sequence of \cref{lemma-ses}. By assumption, ${\rm Der}({\rm Sym}_R(L))$ is of rank $n+d$ over ${\rm Sym}_R(L)$ and this sequence shows that $Q_L\cong\Omega^{n+d}_{{\rm Sym}_R(L)}$ as $R$-modules. 
There is a canonical duality isomorphism 
\[
\mathfrak{X}^\bullet_{{\rm Sym}_R(L)}\otimes_R\Omega^{n+d}_{{\rm Sym}_R(L)}\cong \Omega^{\bullet-n-d}_{{\rm Sym}_R(L)},
\]
which is an algebraic analogue of the duality between Poisson homology and cohomology in \cite{Xu}.
We therefore see that the cochain complex computing the cohomology in \cref{eq-hi} is given by 
\[
\sfOmega_{\rm nl}^{n+d-\bullet}({\rm Sym}({\rm ad})\otimes Q_L)\cong \Omega^\bullet_{{\rm Sym}_R(L)},
\]
with differential equal to $L_P:=d\circ \iota_P+\iota_P\circ d$, with $P\in\mathfrak{X}^2_{{\rm Sym}_R(L)}$ the Poisson bracket \cref{eq-epb} and $d: \Omega^{\bullet}_{{\rm Sym}_R(L)}\to \Omega^{\bullet+1}_{{\rm Sym}_R(L)}$ the K\"ahler differential.

We denote by $(CC_\bullet(\sfU),b,B)$ the standard mixed complex computing Hochschild and 
cyclic homology, c.f. \cite[\S 9.8]{weibel}. The PBW filtration defined in \cref{eq-pbwf} turns the Hochschild complex $(C_\bullet(\sfU),b)$ into a filtered complex and the associated spectral 
sequence has on its first page
\[
E^1_{p,q}=\Omega^{p+q}_{{\rm Sym}_R(L)}(p),
\]
the algebraic $p+q$-forms of polynomial degree $p$. It is well-known that the Poisson bracket \cref{eq-epb} on ${\rm Sym}_R(L)$ is induced by the 
commutator $[D,E]:=DE-ED$ on $\sfU(L,R)$, and therefore the induced differential on the first page is exactly given  by $L_P: \Omega^{p+q}_{{\rm Sym}_R(L)}(p)
\to\Omega^{p+q-1}_{{\rm Sym}_R(L)}(p-1)$, c.f. \cite{BrylinskiP}. (Remark that the Poisson bracket \cref{eq-epb} is of polynomial degree $-1$.)
Since we already know that this cochain complex on the first page computes the Hochschild homology, c.f. \cref{eq-hi}, we conclude that the spectral sequence degenerates 
from this page on. By the SBI sequence relating Hochschild and cyclic homology, the spectral sequence must also degenerate for the cyclic bicomplex. 
Similar considerations as above show that the cyclic differential $B$ induces the K\"ahler differential $d$ on the first page, and this leads to the final conclusion:
\begin{theorem}
\label{th:dht}
Let $(L,R)$ be a Lie--Rinehart algebra with $L$ projective over $R$ of constant rank, and $R$ satisfying van den Bergh duality. 
Then there are canonical isomorphisms
\begin{align*}
HH_\bullet(\sfU(L,R))&\cong H_\bullet(\Omega^\bullet_{{\rm Sym}_R(L)},L_P),\\
HC_\bullet(\sfU(L,R))&\cong H_\bullet(\Omega^\bullet_{{\rm Sym}_R(L)}[u],L_P+ud),
\end{align*} 
with $u$ a formal variable of degree $-2$.
\end{theorem}
\section{Examples} 
We now briefly discuss the implications of our computation for our main examples of Lie--Rinehart algebras: Lie algebras and differential operators, supplemented with two more.
\subsection{Lie algebras} Going back to \cref{ex-la}, we consider the universal enveloping algebra $\sfU(\mathfrak{g})$ of a Lie algebra $\mathfrak{g}$.
In this case $R=k$ and ${\rm Sym}({\rm ad}(\mathfrak{g}))$ is just the adjoint representation extended to the symmetric algebra and \cref{thm-m1} gives
\[
H^\bullet(\sfU(\mathfrak{g}),\sfU(\mathfrak{g}))\cong H^\bullet_{\rm CE}(\mathfrak{g},{\rm Sym}(\mathfrak{g})).
\]
Dually, we use that $Q_\mathfrak{g}={\rm Det}(\mathfrak{g})$ and that the Lie algebra cohomology twisted with $Q_\mathfrak{g}$ identifies with Lie algebra homology: $H^\bullet_{\rm CE}(\mathfrak{g};\mathsf{M}\otimes Q_\mathfrak{g})\cong H^{\rm CE}_{d-\bullet}(\mathfrak{g};\mathsf{M})$. With that, \cref{th:dht} gives the natural isomorphisms
\begin{align*}
HH_\bullet(\sfU(\mathfrak{g}))&\cong H^{\rm CE}_\bullet(\mathfrak{g},{\rm Sym}(\mathfrak{g})),\\
HC_\bullet(\sfU(\mathfrak{g}))&\cong H^{\rm CE}_\bullet(\Omega^\bullet_{{\rm Sym}(\mathfrak{g})}[u],L_P+ud)
\end{align*}
These isomorphisms agree with the results in the literature first obtained in \cite{Kassel}. In \cite{PT}, the isomorphism for Hochschild cohomology is derived 
from Kontsevich' formula for deformation quantization, and shown to be compatible with the natural graded algebra structure on both sides (Poisson and Hochschild cohomology). In this paper we have not considered the algebra structure.

\subsection{Differential operators} For the algebra of differential operators on an affine variety, c.f.\ \cref{ex-do}, our main theorem \cref{thm-m1} states that
the Hochschild cohomology of $\mathcal{D}(\mathsf{X})$ is given by the cohomology of the symmetric powers of the adjoint quasi-module $({\rm Sym}_R({\rm ad}({\rm Der}(R))),-\delta)$ of \cref{section:epc}. But in this case the underlying complex is simply a Koszul complex and therefore the inclusion
\[
(R,0)\hookrightarrow ({\rm Sym}_R({\rm ad}({\rm Der}(R))),-\delta)
\]
a quasi-isomorphism. Therefore we find  that
\[
H^\bullet(\mathcal{D}(\mathsf{X}),\mathcal{D}(\mathsf{X}))\cong H^{\bullet}(\mathsf{X}).
\]
For the homology theory we remark that
\[
Q_{{\rm Der}(R)}:=\bigwedge^n{\rm Der}(R)\otimes\Omega^n_R\cong R,
\]
and therefore \cref{th:dht} gives 
\begin{align*}
HH_\bullet(\mathcal{D}(\mathsf{X}))\cong H^{2n-\bullet}(\mathsf{X}),\\
HC_\bullet(\mathcal{D}(\mathsf{X}))\cong\bigoplus_{k\geq 0}H^{2n+2k-\bullet}(\mathsf{X}).
\end{align*}
These last two isomorphisms were originally obtained in \cite{Wodzicki}
\subsection{Lie algebra actions} Let $\mathfrak{g}$ be a Lie algebra (over the ground field $\mathbb{K}$) that acts on the commutative algebra $R$ by way of a Lie algebra morphism $\rho: \mathfrak{g}\to{\rm Der}(R)$. This leads to a Lie-Rinehart algebra $(L,R)=(R\otimes\mathfrak{g},R)$ with bracket
\begin{equation*}
[r\otimes v,r'\otimes v']=rr'\otimes [v,v']+r\rho(v)(r')\otimes v'-r'\rho(v')(r)\otimes v
\end{equation*}
and anchor $\rho: R\otimes\mathfrak{g}\to\Der(R)$ given by
\begin{equation*}
\rho(r\otimes v)=r\rho(v).
\end{equation*}
Its universal enveloping algebra $\sfU(L,R)$ is a semi-direct product $R\rtimes \sfU(\mathfrak{g})$ where the action of $\sfU(\mathfrak{g})$ on $R$ is given by the map $\sfU(\mathfrak{g})\to \sfU({\rm Der}(R),R)\hookrightarrow{\rm End}(R)$ induced by the action $\rho$.

Now, if we look at the non-linear complex $\sfOmega_{\rm nl}({\rm Sym}({\rm ad}(L)))$, we note that since $R\otimes\mathfrak{g}$ is a free $R$-module, the non-linearity equation \eqref{eq-nlcec} simplifies a lot. Indeed, a pair $(\phi_0,\phi_1)$ of a map
\[\phi_0: R\otimes\mathfrak{g}\to R\otimes\mathfrak{g}\]
with a symbol
\[\phi_1\in \Der(R)\]
such that
\[\phi_0(rr'\otimes v)=r\phi_0(r'\otimes v)+h_{r,r'\otimes v}(\phi_1)\]
is completely characterized by $\phi_1$ and the restriction of $\phi_0$ to $1\otimes\mathfrak{g}$, since
\[\phi_0(r\otimes v)=r\phi_0(1\otimes v)+h_{r,1\otimes v}(\phi_1).\]
The important point is that the restriction of $\phi_0$ and the symbol $\phi_1$ are independent of each other, this restriction gives a bijection
\[\{(\phi_0,\phi_1):\phi_0(rr'\otimes v)=r\phi_0(r'\otimes v)+h_{r,r'\otimes v}(\phi_1)\}\cong {\rm Hom}(\mathfrak{g},R\otimes\mathfrak{g})\oplus{\rm Der}(R).\]

In general, restricting non-linear cochains to $1\otimes\mathfrak{g}$ induces an isomorphism on the level of chains between $\sfOmega_{\rm nl}({\rm Sym}({\rm ad}(L)))$ and $\sfOmega(\mathfrak{g},\bigwedge_R{\rm Der}(R)\otimes{\rm Sym}(\mathfrak{g}))$ where the latter is given by
\[\sfOmega^n(\mathfrak{g},\bigwedge_R{\rm Der}(R)\otimes{\rm Sym}(\mathfrak{g}))=\bigoplus_{i=0}^n {\rm Hom}\left(\bigwedge^{n-i}\mathfrak{g},\bigwedge^i_R{\rm Der}(R)\otimes{\rm Sym}(\mathfrak{g})\right)\]
Since the bracket on $L$ satisfies
\[[1\otimes v,1\otimes v']=1\otimes [v,v']\]
we see that the differential induced on $\sfOmega(\mathfrak{g},\bigwedge_R{\rm Der}(R)\otimes{\rm Sym}(\mathfrak{g}))$ is the one induced by a double complex, with the differential in one direction given by the Chevalley-Eilenberg differential for the $\mathfrak{g}$-representation $\bigwedge_R{\rm Der}(R)\otimes{\rm Sym}(\mathfrak{g})$, and the differential in the other direction is given by postcomposition with the map $\bigwedge^p_R{\rm Der}(R)\otimes{\rm Sym}^q(\mathfrak{g})\to \bigwedge^{p+1}_R{\rm Der}(R)\otimes {\rm Sym}^{q-1}(\mathfrak{g})$ induced by the action $\rho$.

An application of \cref{thm-m1} now gives a spectral sequence
\[
E_1^{p,q}=H^q_{\rm CE}\left(\mathfrak{g},\bigwedge^p_R{\rm Der}(R)\otimes{\rm Sym}(\mathfrak{g})\right)\quad\Longrightarrow\quad HH^{p+q}(R\rtimes\sfU(\mathfrak{g}),R\rtimes\sfU(\mathfrak{g})).
\]
\subsection{Central arrangements of lines}
Finally, we briefly discuss the case of  differential operators tangent to a central arrangement of lines. Its Hochschild cohomology was computed recently in \cite{KSA,KL},
here we mainly point out the Poisson complex that arises from our point of view.

We consider $R=\mathbb{K}[x,y]$ and let $Q=xF$ be the defining polynomial of $r+2$ lines in the plane $\mathbb{A}^2={\rm Spec}(R)$, where we have put the first line at $x=0$.
The Lie subalgebra 
\[
L:=\{X\in{\rm Der}(R),~X(Q)\subset QR\}
\]
forms a Lie-Rinehart algebra which is free over $R$ by Saito's criterium with generators $E:=x\partial_x+y\partial_y$ and $D:=F\partial_y$.
Therefore ${\rm Sym}_R(L)=\mathbb{K}[x,y,D,E]$ and from the elementary Lie brackets between the generators listed in \cite[\S 1]{KSA} we 
infer that the Poisson structure is given by
\begin{equation}
\label{psca}
P:=x\partial_E\wedge\partial_x+F\partial_D\wedge\partial_y+y\partial_E\wedge\partial_y+rD\partial_E\wedge\partial_D.
\end{equation}
With this Poisson structure, the grading of {\em loc. cit.} is given by the Euler vector field $\mathcal{E}:=\{-,E\}$.
On the Poisson complex $(\mathfrak{X}^\bullet_{{\rm Sym}_R(L)},\delta_P:=[-,P])$ we therefore consider the 
operator $s_0:=\iota_{dE}$ given by contraction with $dE\in\Omega^1_{{\rm Sym}_R(L)}$. A short computation shows that
\[
\delta_P\circ s_0+s_0\circ \delta_P=\mathcal{E}.
\]
We can therefore contract the Poisson complex to the degree zero part $(\mathfrak{X}^\bullet_{{\rm Sym}_R(L),0},\delta_P)$. 
This is a small complex in which explicit computations can be made. For example in cochain degree $0$ the only polynomials of degree $0$ are given by $E^q$
for some $q\geq 0$, and we see immediately from \cref{psca} that in $\ker(\delta_P)$ we must have $q=0$ and therefore $HH^0(\sfU,\sfU)=\mathbb{K}$. The
computations for the remaining degrees are very similar to \cite[\S 3]{KSA}. Remark that these authors also consider the Gerstenhaber structure on $HH^\bullet(\sfU,\sfU)$.

\bibliographystyle{alpha}
\bibliography{references}

\end{document}